\documentclass[10pt,a4paper]{article}
\usepackage[T1]{fontenc}
\usepackage[utf8]{inputenc}
\usepackage{authblk}
\usepackage{amsmath,amssymb,amsthm,esint,bm}
\usepackage{mathrsfs}
\usepackage{bookmark}
\usepackage{amsmath}
\usepackage{enumerate}
\allowdisplaybreaks[3]
\allowdisplaybreaks[4]

\newtheorem{definition}{Definition}[section]
\newtheorem{theorem}[definition]{Theorem}
\newtheorem{lemma}[definition]{Lemma}

\newtheorem{corollary}[definition]{Corollary}
\theoremstyle{remark}
\newtheorem{remark}[definition]{Remark}
\numberwithin{equation}{section}

\newcounter{stepnum}

\newcommand{\step}{%
	\par
	\refstepcounter{stepnum}%
	\textbf{Step \arabic{stepnum}}.\enspace\ignorespaces
}
\usepackage{etoolbox}
\AtBeginEnvironment{proof}{\setcounter{stepnum}{0}}

\newcommand{\R}{\mathbb{R}}

\newcommand{\RNum}[1]{\uppercase\expandafter{\romannumeral #1\relax}}

\title{Potential estimates and improved Sobolev regularity for nonlocal fully nonlinear elliptic equations}
\author[a]{Yao Zhang}
\author[a]{Lingwei Ma}
\author[a]{Zhenqiu Zhang\thanks{Corresponding author.}}
\affil[a]{School of Mathematical Sciences and LPMC, Nankai University, Tianjin, 300071, P.R. China}

\date{\today}

\usepackage{hyperref}
\begin{document}
	\maketitle
	\footnotetext[1]{E-mail: 1120250049@mail.nankai.edu.cn(Y. Zhang), malingwei@nankai.edu.cn (L. Ma), zqzhang@nankai.edu.cn (Z. Zhang).}
	
	\maketitle
	\begin{abstract}
		This paper investigates pointwise and oscillatory potential estimates for viscosity solutions and their gradients to a class of nonlocal fully nonlinear elliptic equations. The key ingredient of our argument is to establish precise pointwise bounds for a novel fractional sharp maximal function defined in terms of affine functions. As an immediate consequence, we prove higher order Sobolev regularity results of Calder\'{o}n-Zygmund type for such equations without imposing any global boundedness assumptions on viscosity solutions. The potential theory developed herein extends the work of Daskalopoulos-Kuusi-Mingione (Commun. Partial Differ. Equ. 39: 574-590, 2014) for fully nonlinear second order elliptic equations to the fractional fully nonlinear setting.\\

		
	\end{abstract}
	
	
	\section{Introduction and main results}\label{section1}
	
	This paper is concerned with the pointwise potential estimates for viscosity solutions to the following nonlocal fully nonlinear elliptic  problems
	\begin{equation}\label{prob1}
		\mathcal{I}(u,x)=f(x)\quad \text{in } \Omega,
	\end{equation}
		where $\Omega\subset \mathbb{R}^n\,(n\ge2)$ is a bounded domain and the source term $f\in L^p(\Omega)\cap C(\Omega)$. Considerations arising from two-player zero-sum stochastic games lead to the fully nonlinear integro-differential operator $\mathcal{I}$ defined as
	\begin{equation}\label{i}
		\mathcal{I}(u,x):=\inf_{\beta\in \mathcal{B}}\sup_{\alpha\in\mathcal{A}}\{-\mathcal{L}_{\alpha\beta} u(x)\}
	\end{equation}
	where $\mathcal{A}$, $\mathcal{B}$ are index sets and each $\mathcal{L}_{\alpha\beta}$ is a linear integro-differential elliptic operator of order $2s$ with $s\in(0,1)$, given by
	\begin{equation}\label{k0}
		\mathcal{L}_{\alpha\beta} u(x):=\frac{1}{2} \int_{\mathbb{R}^n}\delta(u,x,y)K_{\alpha\beta}(x,y)dy
	\end{equation}
	with $\delta(u,x,y):=2u(x)-u(x+y)-u(x-y)$ and  $$K_{\alpha\beta}(x,y):=\frac{y^TA_{\alpha\beta}(x)y}{|y|^{n+2s+2}}.$$ Let $S_{\lambda,\Lambda}$ denote the set of real symmetric $n\times n$ matrice $A_{\alpha\beta}$  satisfying the ellipticity condition
	\begin{equation}\label{k2}
	 \lambda I\le A_{\alpha\beta}\le \Lambda I\quad \text{in }\mathbb{R}^n,
	\end{equation}for some $0<\lambda\le \Lambda$.
	For the sake of convenience, we denote by $\mathfrak{L}_s(\lambda,\Lambda)$ and $\mathfrak{J}_s(\lambda,\Lambda)$ the classes of linear and fully nonlinear integro-differential elliptic operators of order 2s, respectively, as follows.
	\begin{definition}
		Let $s\in(0,1)$, $0<\lambda\le \Lambda$ and $\omega$ be a modulus of continuity. Fix $M>0$ and $\eta\in(0,1)$ such that $\eta+2s\notin\mathbb{N}$. Define the classes of linear integro-differential operators associated with the kernerl $K$ and the coefficient matrix $A$ by
		\begin{equation*}
			\mathfrak{L}^0_s(\lambda,\Lambda):=\{\mathcal{L}:\eqref{k0}\text{ and }\eqref{k2}\text{ hold}\}.
		\end{equation*}Define the classes of linear integro-differential  operators with the coefficient matrix $A$ satisfying a modulus of continuity condition by
	\begin{equation*}
	\mathfrak{L}^1_s(\lambda,\Lambda):=\{\mathcal{L}\in\mathfrak{L}^0_s(\lambda,\Lambda) :|A(x)-A(x')|\le \omega(|x-x'|)\text{ for all $x,x'\in \Omega$}\}.
	\end{equation*}Define the classes of linear integro-differential  operators with the coefficient matrix $A$ satisfying $\eta$-H\"{o}lder condition by
\begin{equation*}
\mathfrak{L}^2_s(\lambda,\Lambda):=\{\mathcal{L}\in\mathfrak{L}^0_s(\lambda,\Lambda) :[A]_{C^\eta(\Omega)}\le M\}.
\end{equation*}For $j=0,1,2$, the corresponding fully nonlinear integro-differential  operator classes are defined by
\begin{equation*}
	\mathfrak{J}^j_s(\lambda,\Lambda):=\{\mathcal{I}:\eqref{i}\text{ with $\mathcal{L}\in\mathfrak{L}^j_s(\lambda,\Lambda)$}\}.
\end{equation*}
	\end{definition}


	The considerable interest in nonlocal integro-differential equations is largely driven by their applications in both pure and applied mathematics. Acting as the mathematical realization of anomalous diffusion, L\'{e}vy flights, and long-range interactions \cite{MK, Laskin}, these fractional frameworks naturally arise in modeling discontinuous jump-diffusion stochastic processes \cite{apple}. Furthermore, such nonlocal operators play fundamental role in several physical models \cite{CF00, CA2010, LY,IS22}. In light of this broad range of applications, developing a comprehensive regularity theory for fully nonlinear integro-differential equations is of central importance. The foundational regularity framework for fully nonlinear integro-differential equations \eqref{prob1}
	was established in the celebrated works of Caffarelli and Silvestre \cite{CS2009, CS2011, CS20111}.
	More precisely,the authors of \cite{CS2011} studied integro-differential equations \eqref{prob1} with variable coefficients and rigorously established $C^{1,\alpha}$ regularity under the assumption that the kernel $K(x,y)$  satisfies a modulus of continuity condition in $x$.
  Under an additional concavity assumption on the operator $\mathcal{I}$ with the form
	\begin{equation*}
		\mathcal{I}(u,x):=\inf_{\beta\in \mathcal{B}} \{-\mathcal{L}_{\beta}u(x)\},
	\end{equation*}and the translation-invariant kernel, a nonlocal version of the Evans-Krylov theorem ($C^{2s,\alpha}$-regularity) was obtained in \cite{CS20111}. Later, Serra \cite{serra} strengthened the assumption on the kernel and extended these results to concave equations \eqref{prob1}  that need not be translation-invariant,
	deriving $C^{2s,\alpha}$ estimates. The aforementioned literature has thoroughly addressed the continuity and classical smoothness properties, a natural next step is therefore to investigate the higher-order integrability of the solutions. Nevertheless, the nonlocal Aleksandrov-Bakelman-Pucci (ABP) estimate originally derived in \cite{CS2009} requires the right-hand side to depend on the boundedness of source term $f$. To relax this stringent requirement, Guillen and Schwab \cite{GS2012} generalized the nonlocal ABP estimate to the following form
	\begin{equation}\label{ABP2}
		-\inf_{B_1} u\le \frac{C(n)}{\lambda}\|f\|_{L^\infty(\Gamma_u)}^{1-s}\|f\|_{L^n(\Gamma_u)}^{s},
	\end{equation}where $\Gamma_u\subset B_1$ is a coincidence set with regard to $u$.
	Building upon the extremal estimates \eqref{ABP2}, Yu \cite{yu} subsequently proved the $W^{2s, \epsilon}$ regularity for solutions  which belong to $L^\infty(\mathbb{R}^n)$ to nonlocal elliptic equations for some $\epsilon>0$. More recently, Kitano \cite{Kitano} successfully established nonlocal ABP estimates controlled solely  by the $L^p$ norms of source term $f$, and further obtained the $W^{2s,p}$ regularity results for the essentially bounded solutions.
	
Remarkably, a recent breakthrough by Diening and Nowak \cite{DN2025} reveals a new and elegant approach from the viewpoint of modern harmonic analysis. They showed that Calder\'{o}n-Zygmund type estimates can be derived directly by combining pointwise sharp maximal function estimates with  potential bounds. Motivated by this approach, the main objective of this paper is to establish potential estimates for viscosity solutions and their gradients in the setting of nonlocal fully nonlinear equations, and thereby derive Calder\'{o}n-Zygmund type estimates in a natural way.

	The study of pointwise potential estimates has undergone profound development in the regularity theory of partial differential equations, deeply rooted in the techniques of harmonic analysis. For local elliptic equations, potential theory is well-established for both divergence and non-divergence form operators. The fundamental results due to  Kilpel\"{a}inen and Malý \cite{KM}, who investigated a class of divergence-form quasilinear second-order elliptic equations and obtained pointwise zero-order potential estimates. Further remarkable results for the gradient of solutions have been achieved by  Mingione \cite{Min11}. In a vast body of subsequent literature, similar potential estimates have been successfully extended to more general local nonlinear second-order elliptic and parabolic equations, refer to \cite{DM10, KM2012, KM14, KM18, BY19, Ma1, Ma2} and the reference therein. For local fully nonlinear equations, Daskalopoulos, Kuusi and Mingione \cite{DKM2014} achieved a breakthrough by establishing gradient potential estimates and simultaneously obtaining the Hölder continuity of the gradient for the viscosity solutions to
	\begin{equation*}
		F(D^2u,x)=f(x)\quad\text{in $\Omega$},
	\end{equation*}where $f\in L^p(\Omega)$ with $p>n$. Subsequently, Pimentel and Walker \cite{PW2022} generalized their results to more general fully nonlinear elliptic equations with bounded measurable ingredients
	\begin{equation*}
		F(D^2u,Du,u,x)=f\quad\text{in }\Omega,
	\end{equation*}with bounded measurable ingredients. Nonetheless, their approach relies heavily on the established $W^{1,q}$ estimates of the $L^p$-viscosity solutions. For further results on the regularity theory of local fully nonlinear elliptic equations, we refer to \cite{1989,user,caffarelli1995fully,1996,winter,ZY}.
	
	Regarding nonlocal elliptic equations, Kuusi, Mingione and Sire \cite{KMS2015} investigated a class of fractional $p$-Laplace equations and introduced the nonlocal excess functional and obtained zero-order potential estimates. Due to the lack of  $C^{1,\alpha }$ estimates for solutions to homogeneous nonlocal equations, establishing gradient potential estimates for the solutions encounters substantial additional difficulties. Later, Kuusi, Nowak and Sire \cite{KNS2022} studied a class of nonlocal linear elliptic equations with kernels satisfying suitable regularity assumptions and successfully derived gradient potential estimates. In \cite{D2024}, the authors achieved gradient potential estimates for nonlinear nonlocal equations under Lipschitz and monotonicity assumptions. More recently, building upon the $C^{1,\alpha}$ regularity for solutions to the homogeneous fractional $p$-Laplace equation established in \cite{1}, Li and Zhang \cite{2} derive pointwise gradient potential estimates for the fractional $p$-Laplace elliptic equation with measure data in the regime $p>2$ and $s\in(\frac{p-1}{p},1)$. Furthermore, pointwise and universal potential estimates for solutions and their gradients to mixed local and nonlocal elliptic and parabolic equations can be found in \cite{3,4,5}. However, all the aforementioned methodologies rely fundamentally on the variational structure of the equations, rendering them inapplicable to non-divergence fully nonlinear elliptic equations. Furthermore, unlike the local case, the requisite differentiability estimates for possibly unbounded viscosity solutions of fully nonlinear nonlocal equations remain unavailable, posing a major theoretical barrier.

To the best of our knowledge, there is currently no existing work concerning the potential estimates for viscosity solutions to fully nonlinear integro-differential equations. This paper aims to bridge this theoretical gap. Adapting and profoundly extending the potential estimation techniques from local contexts, we design a refined Campanato iterative approximation scheme tailored to handle the singular nonlocal tails and long-range interactions. We successfully establish pointwise potential bounds for both the viscosity solution $u$ and its gradient $Du$, as well as pointwise estimates for suitably defined nonlocal fractional sharp maximal functions associated with them. Having secured the bounds for the fractional sharp maximal functions, we further employ harmonic analysis tools to ultimately achieve the pointwise $C^{1,\alpha}$ continuity and the Calder\'{o}n-Zygmund estimates of the viscosity solutions.
	
	Before stating our main results, we introduce several essential notations and function spaces that will be used throughout the paper. Due to the nonlocal nature of the integro-differential operators, functions are required to exhibit a prescribed growth control at infinity to ensure that the singular integrals are well-defined. To this end, we work within the weighted Lebesgue space $L^1_{2s}(\mathbb{R}^n)$, which consists of all measurable functions $u$ satisfying
	\begin{equation}
		\|u\|_{L^1_{2s}(\mathbb{R}^n)} := \int_{\mathbb{R}^n} \frac{|u(y)|}{1+|y|^{n+2s}} dy < +\infty.
	\end{equation}
	Associated with this function space, we define the nonlocal tail term. For a function $u \in L^1_{2s}(\mathbb{R}^n)$, its tail captures the far field contributions outside a local ball $B_r(x_0)\Subset \Omega$, and is given by
	\begin{equation}
		\operatorname{Tail}(u;x_0,r) := r^{2s} \int_{\mathbb{R}^n\setminus B_r(x_0)} \frac{|u(y)|}{|x_0-y|^{n+2s}} dy.
	\end{equation}
	It is worth noting that the assumption $u\in L_{2s}^1(\R^n)$ naturally guarantees the finiteness of the tail term for any $r>0$.
	
	For an $L^p$ integrable function $f$, the $L^p$ version of the truncated Riesz potential of $f$ with order $\beta>0$ is defined as follows
	\begin{equation}
		\mathbf{I}^f_{\beta,p}(x_0,r) := \int^r_0 \rho^{\beta} \left( \fint_{B_\rho(x_0)} |f(y)|^p \, dy \right)^{1/p} \frac{d\rho}{\rho}.
	\end{equation}

	We first present our result on the pointwise potential estimate for the solution itself.
	\begin{theorem}\label{thm1.1}
		Let $s\in(0,1)$, $p^\ast=p^\ast(n,\lambda,\Lambda)\in [n/2,n)$ and $f\in L^p(\Omega)\cap C(\Omega)$ with $p>\max\{p^{\ast},n/2s\}$. Assume that $u\in C(\Omega)\cap L^1_{2s}(\mathbb{R}^n)$ is a viscosity solution to \eqref{prob1} with the form \eqref{i} and $\mathcal{I}\in \mathfrak{J}^0_s(\lambda,\Lambda)$. Then there exists a constant $C=C(n,s,p,\lambda,\Lambda)>0$ such that for every $B_R(x)\Subset \Omega$, it holds that
		\begin{equation*}
			|u(x)|\le C\left(\mathbf{I}^{f}_{2s,p}(x,R)+\fint_{B_R(x)}|u(y)|dy+\operatorname{Tail}(u;x,R)\right).
		\end{equation*}
	\end{theorem}

	Building upon pointwise estimates for the appropriate sharp maximal operators of the solution, we derive the universal potential estimate of the solution. Indeed, the result below establishes a connection between Krylov-Safonov theory and nonlinear potential theory.
	\begin{theorem}\label{thm1.3}
		Let $s\in(0,1)$, $p^\ast=p^\ast(n,\lambda,\Lambda)\in [n/2,n)$,  $B_R(x_0)\Subset \Omega$ and $f\in L^p(\Omega)\cap C(\Omega)$ with $p>\max\{p^{\ast},n/2s\}$. Assume that $u\in C(\Omega)\cap L^1_{2s}(\mathbb{R}^n)$ is a viscosity solution to \eqref{prob1} with the form \eqref{i} and $\mathcal{I}\in \mathfrak{J}^0_s(\lambda,\Lambda)$. Then for any $\tilde{\alpha}\in (0,\alpha_0)$, there exists a constant $C=C(n,s,p,\lambda,\Lambda,\tilde{\alpha})>0$ such that
		\begin{equation*}
			|u(x)-u(y)|\le C\left[\mathbf{I}^{f}_{2s-\alpha,p}(x,R)+\mathbf{I}^{f}_{2s-\alpha,p}(y,R)+\frac{1}{R^\alpha}\left(\fint_{B_R(x_0)}|u|dx+\operatorname{Tail}(u;x_0,R)\right)\right]|x-y|^\alpha,
		\end{equation*}holds uniformly in $0\le\alpha\le\tilde{\alpha}$ for every $x,y\in B_{R/8}(x_0)$, where $\alpha_0=\alpha_0(n,s,\lambda,\Lambda)\in (0,\min\{2s,1\})$ is the H\"{o}lder continuity exponent for solutions to the corresponding homogeneous equation $\mathcal{I}(v,x)=0$, given by Lemma \ref{holder}.
	\end{theorem}

	Subsequently, for $s\in(1/2,1)$, by imposing additional regularity assumptions on the coefficient of the equation, we further obtain pointwise and oscillation-type potential estimates for the gradient of the solution, respectively.
	\begin{theorem}\label{thm1.4}
		
		Let $s\in (1/2,1)$ and $f\in L^p(\Omega)\cap C(\Omega)$ with $p>n/(2s-1)$. Assume that $u\in C(\Omega)\cap L^1_{2s}(\mathbb{R}^n)$ is a viscosity solution to \eqref{prob1} with the form \eqref{i} and $\mathcal{I}\in \mathfrak{J}^1_s(\lambda,\Lambda)$. Then for every $B_R(x)\Subset \Omega$, u is Fr\'{e}chet differentiable at $x$. Moreover, there exists a constant $C=C(n,s,p,\lambda,\Lambda,\omega(\cdot))>0$ holds that
		\begin{equation*}
			|Du(x)|\le C\mathbf{I}^{f}_{2s-1,p}(x,R)+CR^{-1}\left(\fint_{B_{R}(x)}|u(y)|dy+\operatorname{Tail}(u;x,R)\right).
		\end{equation*}
	\end{theorem}
	\begin{theorem}\label{thm1.5}
		Let $s\in (1/2,1)$, $B_R(x_0)\Subset\Omega$ and $f\in L^p(\Omega)\cap C(\Omega)$ with $p>n/(2s-1)$. Assume that $u\in C(\Omega)\cap L^1_{2s}(\mathbb{R}^n)$ is a viscosity solution to \eqref{prob1} with the form \eqref{i} and $\mathcal{I}\in \mathfrak{J}^1_s(\lambda,\Lambda)$. Then for any $\tilde{\alpha}\in (0,\min\{\bar{\alpha},2s-1\})$, there exists a constant $C=C(n,s,p,\lambda,\Lambda,\tilde{\alpha},\omega(\cdot))>0$ such that
		\begin{equation*}
			\begin{aligned}
				|Du(x)&-Du(y)|\\
				&\le C\left[\mathbf{I}^{f}_{2s-1-\alpha,p}(x,R)+\mathbf{I}^{f}_{2s-1-\alpha,p}(y,R)+\frac{1}{R^{1+\alpha}}\left(\fint_{B_{R}(x_0)}|u|dy+\operatorname{Tail}(u;x_0,R)\right)\right]|x-y|^\alpha,
			\end{aligned}
		\end{equation*}holds uniformly in $0\le\alpha\le\tilde{\alpha}$ for every $x,y\in B_{R/8}(x_0)$, where $\bar{\alpha}=\bar{\alpha}(n,s,\lambda,\Lambda)\in (0,\alpha_0)$ denotes the H\"{o}lder continuity exponent for gradients of solutions to the corresponding homogeneous equation $\mathcal{I}(v,x)=0$, given by Lemma \ref{1+alpha}.
	\end{theorem}

As a by-product of the gradient potential estimate from Theorem \ref{thm1.4}, and the pointwise estimate for the certain sharp maximal operators of the gradient, fine higher regulairy results of Calder\'{o}n-Zygmund type can be proved for fully nonlinear nonlocal elliptic equations.
	\begin{corollary}\label{cor1.7}
		Let $s\in (1/2,1)$. Assume that $u\in C(\Omega)\cap L^1_{2s}(\mathbb{R}^n)$ is a viscosity solution to \eqref{prob1} with the form \eqref{i} and $\mathcal{I}\in \mathfrak{J}^1_s(\lambda,\Lambda)$. For every $t\in(0,\min\{\bar{\alpha},2s-1\})$ assume that $f\in L^q(\Omega)\cap C(\Omega)$ with $n/(2s-1)< q<n/(2s-1-t)$ and denote $\bar{q}:=\frac{nq}{n-(2s-1-t)q}$. Then $u\in W^{1+t,\bar{q}}_{loc}(\Omega)$ and for $B_R(x_0)\Subset\Omega$ there exists a positive constant $C=C(n,s,q,\lambda,\Lambda,\omega(\cdot),t)>0$, it holds that
			\begin{equation*}
				\begin{aligned}
					&R^{-t}\|Du\|_{L^{\bar{q}}(B_{R/16}(x_0))}+[Du]_{W^{t,\bar{q}}(B_{R/16}(x_0))}\\
					\le& CR^{\frac{n}{\bar{q}}-t-1}\left(\fint_{B_R(x_0)}|u|dy+\operatorname{Tail}(u;x_0,R)\right)+C\|f\|_{L^{q}(B_{R}(x_0))},
				\end{aligned}
			\end{equation*}where $\bar{\alpha}=\bar{\alpha}(n,s,\lambda,\Lambda)\in (0,\alpha_0)$ is given by Lemma \ref{1+alpha}.
	\end{corollary}
\begin{remark}
	All the results stated above concern viscosity solutions of fully nonlinear integro-differential equations with the form \eqref{i}, whose operators are generally neither convex nor concave. In contrast to the results of Daskalopoulos, Kuusi and Mingione \cite{DKM2014} for local fully nonlinear equations, we provide the gradient potential estimates in the nonlocal setting. For the variable-coefficient equation \eqref{prob1}, building on the $C^{1,\alpha}$ regularity result of Caffarelli and Silvestre \cite{CS2011} for bounded source terms, we establish the oscillation estimates of gradients of viscosity solutions with possibly unbounded right-hand sides. While Diening and Nowak \cite{DN2025} established gradient potential estimates and higher differentiability in the nonlocal variational setting, our results address fully nonlinear equations that lack the corresponding energy structure and yield higher differentiability of the viscosity solutions.
	
	If instead, we consider the more structured nonlocal operators
	\begin{equation*}
		\mathcal{I}(u,x):=\inf_{\beta\in \mathcal{B}} \{-\mathcal{L}_{\beta}u(x)\}\quad\text{ with $\mathcal{L}_\beta\in\mathfrak{L}^2_s(\lambda,\Lambda)$},
	\end{equation*}
	then the operator enjoys a concavity structure.
	Using the $C^{2s,\alpha}$ regularity result in Serra \cite{serra}, established under suitable H\"{o}lder continuity assumptions on the data, we exploit the additional concavity structure to obtain improved regularity for equations with possibly unbounded right-hand sides. Specifically, the admissible H\"{o}lder continuity exponents for the gradient of the solutions and the Sobolev regularity ranges in the Calder\'{o}n-Zygmund estimates for $u\in W_{\mathrm{loc}}^{1+t,\bar q}$ can be extended to any $1+t<2s$.
\end{remark}

The remainder of this paper is organized as follows. Section \ref{section2} introduces the basic notations and collects preliminary results. Section \ref{section3} establishes the comparison estimates and decay properties of the global excess functional and we derive pointwise and oscillation potential estimates for viscosity solutions for nonlocal fully nonlinear equation \eqref{prob1} with the form \eqref{i}. Section \ref{section4} develops the analogous theory for gradients, together with the corresponding Calder\'{o}n-Zygmund estimates, based on a global excess functional defined through affine function. Finally, Section \ref{S6} addresses concave nonlocal fully nonlinear equation \eqref{prob1} with the form \eqref{ii} and establishes improved gradient regularity and Calder\'{o}n-Zygmund estimates.

	\section{Preliminaries}\label{section2}
	In this section, we collect basic definitions and technical lemmas required for the proofs of our main theorems in the sequel.
	\subsection{Viscosity solutions and preliminary regularity results}

	We first present the definition of $C$-viscosity subsolution or supersolution to the nonlocal fully nonlinear elliptic equation \eqref{prob1}, and refer the reader to \cite{CS2009,Otonbook} for details.
	\begin{definition}
		Let $s\in(0,1)$, the integro-differential operator $\mathcal{I}\in\mathfrak{J}^0_s(\lambda,\Lambda)$ and the nonhomogeneous term $f\in C(\Omega)$.
		\begin{enumerate}[(i)]
			\item A function $u\in USC(\Omega)\cap L^1_{2s}(\mathbb{R}^n)$ is said to be a viscosity subsolution to \eqref{prob1}, if for any $x\in\Omega$ and any neighborhood $N_x$ of $x$ in $\Omega$, and for every $\phi\in L^1_{2s}(\mathbb{R}^n)$ such that $\phi\in C^2(\overline{N_x})$, $\phi(x)=u(x)$, $\phi\ge u$ in $N_x$ and $\phi= u$ in $\R^n\setminus N_x$, then $\mathcal{I}(\phi,x)\ge f(x)$.
			\item A function $u\in LSC(\Omega)\cap L^1_{2s}(\mathbb{R}^n)$ is said to be a viscosity supersolution to \eqref{prob1}, if for any $x\in\Omega$ and any neighborhood $N_x$ of $x$ in $\Omega$, and for every $\phi\in L^1_{2s}(\mathbb{R}^n)$ such that $\phi\in C^2(\overline{N_x})$, $\phi(x)=u(x)$, $\phi\le u$ in $N_x$ and $\phi= u$ in $\R^n\setminus N_x$, then $\mathcal{I}(\phi,x)\le f(x)$.
			\item If  $u\in C(\Omega)\cap L^1_{2s}(\mathbb{R}^n)$  is simultaneously a viscosity subsolution and a supersolution to \eqref{prob1}, it is called that  $u$  is a viscosity solution to \eqref{prob1}.
		\end{enumerate}
	
	\end{definition}

	The next result shows the relationship between the fully nonlinear integro-differential operator and the nonlocal extremal operator. The fractional Pucci maximal and minimal operator is defined by
	\begin{equation}\label{M+}
		\mathcal{P}^{+}_{\mathfrak{L}^0_s(\lambda,\Lambda)}u(x):=\sup_{\mathcal{L}\in\mathfrak{L}^0_s(\lambda,\Lambda) }\{-\mathcal{L}u\}=\sup_{A\in S_{\lambda,\Lambda}}-\frac{1}{2}\int_{\mathbb{R}^n}\delta(u,x,y)\frac{y^TAy}{|y|^{n+2s+2}}dy
	\end{equation}
	and
	\begin{equation}\label{M-}
		\mathcal{P}^{-}_{\mathfrak{L}^0_s(\lambda,\Lambda)}u(x):=\inf_{\mathcal{L}\in\mathfrak{L}^0_s(\lambda,\Lambda) }\{-\mathcal{L}u\}=\inf_{A\in S_{\lambda,\Lambda}}-\frac{1}{2}\int_{\mathbb{R}^n}\delta(u,x,y)\frac{y^TAy}{|y|^{n+2s+2}}dy,
	\end{equation}respectively. Throughout our article, we denote $\mathcal{P}^{+}:=\mathcal{P}^{+}_{\mathfrak{L}^0_s(\lambda,\Lambda)}$ and $\mathcal{P}^{-}:=\mathcal{P}^{-}_{\mathfrak{L}^0_s(\lambda,\Lambda)}$ for convenience.
	\begin{lemma}\cite[Proposition 3.2.15]{Otonbook}\label{elliptic}
		Let $s\in (0,1)$, $\Omega\subset \mathbb{R}^n$ be any open set and $\mathcal{I}\in\mathfrak{J}^0_s(\lambda,\Lambda)$. Assume that $u,v\in C(\Omega)\cap L^1_{2s}(\mathbb{R}^n)$ satisfying
		\begin{equation*}
			\mathcal{I}(u,x)= f(x)\quad \text{and}\quad \mathcal{I}(v,x)= g(x)\quad\text{in }\Omega
		\end{equation*}in the viscosity sense with  $f,g\in C(\Omega)$. Then
		\begin{equation*}
			\mathcal{P}^{-}(u-v)\le f-g\le \mathcal{P}^{+}(u-v)\quad\text{in }\Omega
		\end{equation*}in the viscosity sense.
	\end{lemma}
	A cornerstone in the regularity theory of non-divergence form equations with bounded measurable coefficients is the local boundedness.
	\begin{lemma}\cite[Theorem 3.3.5]{Otonbook}\label{sup}
		Let $s\in(0,1)$. Assume that $u\in USC(B_1)\cap L^1_{2s}(\mathbb{R}^n)$ satisfies
		\begin{equation*}
			\mathcal{P}^{+}u\ge-C_0\quad\text{in }B_1
		\end{equation*}in the viscosity sense for some $C_0\ge 0$. Then there exists a constant $C=C(n,s,\lambda,\Lambda)$ such that
		\begin{equation*}
			\sup_{B_{1/2}}u\le C\left(\|u\|_{L^1_{2s}(\mathbb{R}^n)}+C_0\right).
		\end{equation*}
	\end{lemma}
	Building directly upon the half Harnack inequality, we state the  H\"{o}lder continuity for viscosity solutions to the extremal equations with bounded measurable coefficients.
	\begin{lemma}\cite[Theorem 3.3.6]{Otonbook}\label{holder}
		Let $s\in(0,1)$. Assume that $u\in C(B_1)\cap L^1_{2s}(\mathbb{R}^n)$ is any viscosity solution to the equation
		\begin{equation*}
			\left\{
			\begin{array}{ll}
				\mathcal{P}^{+}u\ge-C_0 & \text{in } B_1, \\
				\mathcal{P}^{-}u\le C_0 &\text{in } B_1,
			\end{array}
			\right.
		\end{equation*}for some $C_0\ge 0$. Then	$u\in C^{\alpha_0}_{loc}(B_1)$ with
		\begin{equation*}
			\|u\|_{C^{\alpha_0}(B_{1/2})}\le C\left(\|u\|_{L^1_{2s}(\mathbb{R}^n)}+C_0\right),
		\end{equation*}where $C=C(n,s,\lambda,\Lambda)>0$ and $\alpha_0=\alpha_0(n,s,\lambda,\Lambda)\in (0, \min\{2s,1\})$.
	\end{lemma}
In addition, we need the following generalized version of the known nonlocal ABP estimate, which relaxes the global boundedness assumption on solutions from \cite[Theorem 3.1]{Kitano}.
	\begin{lemma}\label{ABP}
		Let $s\in(0,1)$, $p^\ast=p^\ast(n,\lambda,\Lambda)\in [n/2,n)$ and $p>\max\{p^{\ast},n/2s\}$.
		\begin{enumerate}[(i)]
			\item Assume that $u\in LSC(\overline{B_1})\cap L^1_{2s}(\mathbb{R}^n)$ is a viscosity supersolution of
			\begin{equation*}
				\left\{
				\begin{array}{ll}
					\mathcal{P}^{-}u\le f & \text{in } B_1, \\
					u\ge 0 &\text{in } \mathbb{R}^n\setminus B_1,
				\end{array}
				\right.
			\end{equation*}with $f\in C(\overline{B_1})$. Then there exists a constant $C=C(n,p,s,\lambda,\Lambda)$ such that
			\begin{equation*}
				-\inf_{B_1}u\le C\|f^{+}\|_{L^p(B_1)},
			\end{equation*}where $f^+:=\max\{f,0\}$ is the positive part of $f$.
			\item Assume that $u\in USC(\overline{B_1})\cap L^1_{2s}(\mathbb{R}^n)$ is a viscosity supersolution of
			\begin{equation*}
				\left\{
				\begin{array}{ll}
					\mathcal{P}^{+}u\ge f & \text{in } B_1, \\
					u\le 0 &\text{in } \mathbb{R}^n\setminus B_1,
				\end{array}
				\right.
			\end{equation*}with $f\in C(\overline{B_1})$. Then there exists a constant $C=C(n,p,s,\lambda,\Lambda)$ such that
			\begin{equation*}
				\sup_{B_1}u\le C\|f^{-}\|_{L^p(B_1)},
			\end{equation*}where $f^-:=\max\{-f,0\}$ is the negative part of $f$.
		\end{enumerate}
	\end{lemma}
	\begin{proof}
		For any integer $k>0$ and $x\in\mathbb{R}^n$, we define the truncated function
		\begin{equation*}
			u_k(x):=\min\{u(x),k\}.
		\end{equation*}Since $u\in LSC(\overline{B_1})$ and $u\ge 0$ outside $B_1$, we ensure that $u_k\in L^\infty(\mathbb{R}^n)$ and $u_k\ge 0$ in $\mathbb{R}^n\setminus B_1$.
	
		We claim that $u_k\in LSC(\overline{B_1})\cap L^\infty(\mathbb{R}^n) $ is a viscosity supersolution to $\mathcal{P}^- u_k \le f^+$ in $B_1$. Let $x_0 \in B_1$, $N_{x_0}\subset B_1$ be any neighborhood of $x_0$ and for any $\phi\in L^1_{2s}(\mathbb{R}^n)$ such that $\phi\in C^2(\overline{N_{x_0}})$, $\phi(x_0)=u_k(x_0)$ and $\phi(x)\le u_k(x)\text{ in }N_{x_0}.$
		We denote
		\begin{equation*}
			\Phi_k(x):=\left\{
			\begin{array}{ll}
				\phi(x) & x\in N_{x_0}, \\
				u_k(x) &x\notin N_{x_0},
			\end{array}
			\right.
			\quad\text{and}\quad\Phi(x):=\left\{
			\begin{array}{ll}
				\phi(x) & x\in N_{x_0}, \\
				u(x) &x\notin N_{x_0}.
			\end{array}
			\right.
		\end{equation*}If $u_k(x_0)=k$, according to $\phi\le u_k\le k$ in $N_{x_0}$, we obtain $\Phi_k(x)\le k=\Phi_k(x_0)$ for any $x\in \mathbb{R}^n$. Then $\delta(\Phi_k,x_0,y)\ge 0$ and further we have $\mathcal{P}^{-}\Phi_k(x_0)\le 0\le f^+(x_0)$. While if $u_k(x_0)<k$, since $\phi\le u_k\le u$ in $N_{x_0}$ and $\Phi\in L^1_{2s}(\R^n)$, we have $\mathcal{P}^-\Phi(x_0)\le f(x_0)$ according to the definition of viscosity solution. By $\Phi_k(x_0)=\Phi(x_0)$ and $\Phi_k\le \Phi$ in $\mathbb{R}^n$, it is easy to know that
		\begin{equation*}
			\delta(\Phi_k,x_0,y)\ge\delta(\Phi,x_0,y)\quad\text{for any $y\in \mathbb{R}^n$}.
		\end{equation*}Then we obtain $\mathcal{P}^-\Phi_k(x_0)\le\mathcal{P}^-\Phi(x_0) \le f(x_0)\le f^+(x_0)$. Combining the two cases above, we conclude that the claim is valid.
	
		Now, $u_k$ satisfies all the hypotheses of the $L^\infty$-version of the ABP estimate, see \cite[Theorem 3.1]{Kitano}. Then it follows that
		\begin{equation*}
			-\inf_{B_1}u_k\le C\|f^{+}\|_{L^p(B_1)},
		\end{equation*}where the constant $C$ depends only on $n,s,p,\lambda,\Lambda$ but not on $k$. Since we are evaluating the negative infimum and $u_k\le u$, we obtain $-\inf\limits_{B_1} u \le  -\inf\limits_{B_1} u_k$ and finish the first statement. Finally, by replacing $u$ by $-u$, we immediately conclude that $
		\sup\limits_{B_1}u\le C\|f^{-}\|_{L^p(B_1)}$.\end{proof}
To investigate the gradient potential estimates, we need to consider more regular integro-differential operators belonging to the class $\mathfrak{J}^1_s(\lambda,\Lambda)$,  and thereby establish the H\"{o}lder continuity of the gradients of viscosity solutions to homogeneous nonlocal fully nonlinear equations.
	\begin{lemma}\label{1+alpha}
		Given $1/2<s_0\le s<1$, let $\alpha_0>0$ be given by Lemma \ref{holder}. Assume that $u\in C(B_1)\cap L^1_{2s}(\mathbb{R}^n)$ is a viscosity solution of $\mathcal{I}(u,x)=0$ in $B_1$ with the form \eqref{i} and $\mathcal{I}\in \mathfrak{J}^1_s(\lambda,\Lambda)$. Then for some small $\bar{\alpha}\in (0,\alpha_0)$, $u\in C^{1,\bar{\alpha}}_{loc}(B_1)$ and
		\begin{equation*}
			\|u\|_{C^{1,\bar{\alpha}}(B_{1/2})}\le C\|u\|_{L^1_{2s}(\mathbb{R}^n)},
		\end{equation*} where the positive constant $C=C(n,s,\lambda,\Lambda,s_0,\omega(\cdot))$.
	\end{lemma}
	\begin{proof}
		We first consider the cutoff function $0\le\psi\le 1$ such that
		\begin{equation*}
			\psi=1\text{ in $B_{3/4}$}\quad\text{and}\quad\operatorname{supp}\psi\Subset B_{4/5}.
		\end{equation*}Denote $u_1:=\psi u$ and $u_2:=u-u_1$, by using Lemma \ref{elliptic}, then it holds that
		\begin{equation}\label{Pu}
			\mathcal{P}^{-}u\le 0\le  \mathcal{P}^+u \quad\text{in $B_1$}
		\end{equation}in the viscosity sense. Applying Lemma \ref{sup} to \eqref{Pu}, we derive
	\begin{equation*}
	\|u_1\|_{L^\infty(\mathbb{R}^n)}\le \|u\|_{L^\infty(B_{4/5})}\le C\|u\|_{L^1_{2s}(\mathbb{R}^n)}
\end{equation*}from a covering argument. Next we denote
\begin{equation*}
	\tilde{\mathcal{I}}(u_1,x):=\inf_{\beta\in \mathcal{B}}\sup_{\alpha\in\mathcal{A}}\{-\mathcal{L}_{\alpha\beta}u_1(x)+c_{\alpha\beta}(x)\}=\mathcal{I}(u,x)=0\quad\text{in $B_{2/3}$}
\end{equation*}in the viscosity sense, where $c_{\alpha\beta}(x):=-\mathcal{L}_{\alpha\beta}u_2(x)$. Since $u_2=0$ in $B_{3/4}$, by a change of variables,  we deduce that
	\begin{equation*}
		|c_{\alpha\beta}(x)|\le C\|u_2\|_{L^1_{2s}(\mathbb{R}^n)}\le C\|u\|_{L^1_{2s}(\mathbb{R}^n)}
	\end{equation*}for any $\alpha\in\mathcal{A}$, $\beta\in\mathcal{B}$. Thus, we use \cite[Theorem 63]{CS2011} with $u_1$ to obtain that
\begin{equation*}
	\|u\|_{C^{1,\bar{\alpha}}(B_{1/2})}=\|u_1\|_{C^{1,\bar{\alpha}}(B_{1/2})}\le C\left(\|u_1\|_{L^\infty(\mathbb{R}^n)}+\|u\|_{L^1_{2s}(\mathbb{R}^n)}\right)\le C\|u\|_{L^1_{2s}(\mathbb{R}^n)}.
\end{equation*}We finally prove Lemma \ref{1+alpha}.
	\end{proof}

	\subsection{Fractional Sobolev spaces and maximal operators}
	As a direct application of pointwise potential estimates, we prove regularity of viscosity solutions to nonlocal fully nonlinear elliptic equations in suitable fractional Sobolev spaces. For this purpose, we first introduce the definition of the fractional Sobolev spaces $W^{s,p}$.
	\begin{definition}
		Let $\Omega\subset \mathbb{R}^n$ be a domain. For $p\in[1,\infty)$ and $s\in(0,1)$, the fractional Sobolev space is defined as
		\begin{equation*}
			W^{s,p}(\Omega):=\{u\in L^p(\Omega):[u]_{W^{s,p}(\Omega)}<\infty\}
		\end{equation*}with norm
	\begin{equation*}
		\|u\|_{W^{s,p}(\Omega)}:=\|u\|_{L^p(\Omega)}+[u]_{W^{s,p}(\Omega)}.
	\end{equation*}Here
	\begin{equation*}
		[u]_{W^{s,p}(\Omega)}:=\left(\int_{\Omega}\int_{\Omega}\frac{|u(x)-u(y)|^p}{|x-y|^{n+sp}}dydx\right)^{1/p}.
	\end{equation*}
	\end{definition}
	Next, we present the definitions of several fractional maximal operators, which serve as key tools for establishing universal potential estimates.
	\begin{definition}
		Let $x_0\in \mathbb{R}^n$ and $R>0$. For $g\in L^p(B_R(x_0))$, the truncated fractional maximal function of order $\beta\in [0,n]$ is defined as
		\begin{equation*}
			\mathbf{M}_{\beta}[g](x_0,R):=\sup_{0<r\le R}r^{\beta}\fint_{B_r(x_0)}|g(y)|dy,
		\end{equation*}and the $L^p$ version of the truncated fractional maximal function of order $\beta\in [0,n/p]$ is defined as
		\begin{equation*}
			\mathbf{M}_{\beta,p}[g](x_0,R):=\sup_{0<r\le R}r^{\beta}\left(\fint_{B_r(x_0)}|g(y)|^pdy\right)^{\frac{1}{p}}.
		\end{equation*}Moreover, the truncated fractional sharp maximal function of order $\beta\in [0,1]$ is defined as
		\begin{equation*}
			\mathbf{M}_{\beta}^{\#}[g](x_0,R):=\sup_{0<r\le R}r^{-\beta}\fint_{B_r(x_0)}|g(y)-(g)_{B_r(x_0)}|dy,
		\end{equation*}
	where  $(g)_{B_r(x_0)} := \fint_{B_r(x_0)} g(y) dy$ denotes the integral average of $g$ over $B_r(x_0)$.
		When $\beta=0$, we denote $\mathbf{M}_{0}:=\mathbf{M}$ and $\mathbf{M}_{0,p}:=\mathbf{M}_p$ for convenience.
	\end{definition}
Now we establish a connection between the $L^p$ version of the truncated fractional maximal operator and that of the Riesz potential.
	\begin{lemma}\label{Riesz}
		Given $f\in L^p(\Omega)$, $\sigma\in(0,1)$, $R>0$ and $\beta\in[0,n/p]$. For any $B_R(x)\in\Omega$, there exists a positive constant $C=C(n,p,\sigma,\beta)>0$ which holds that
		\begin{equation*}
			\mathbf{M}_{\beta,p}[f](x,\sigma R)\le C\mathbf{I}^f_{\beta,p}(x,R).
		\end{equation*}
	\end{lemma}

	\begin{proof}
		For any $B_R(x)\in \Omega$ and $r\in(0,\sigma R]$, we have $(r,r/\sigma)\subset (0,R]$. Then
		\begin{equation*}
			\begin{aligned}
				\mathbf{I}^f_{\beta,p}(x,R)&=\int_{0}^{R}\rho^{\beta-1}\left(\fint_{B_\rho(x)}|f(y)|^pdy\right)^{1/p}d\rho\\
				&\ge C\int_{r}^{r/\sigma}\rho^{\beta-1-\frac{n}{p}}\|f\|_{L^p(B_\rho(x))}d\rho\\
				&\ge Cr^\beta\left(\fint_{B_r(x)}|f(y)|^pdy\right)^{1/p},
			\end{aligned}
		\end{equation*}where $C=C(n,p,\sigma,\beta)>0$. Since $r\in (0,\sigma R]$ is arbitrary, we obtain that $\mathbf{M}_{\beta,p}[f](x,\sigma R)\le C\mathbf{I}^f_{\beta,p}(x,R)$.\end{proof}

	We further derive the mapping properties of the $L^p$ version of the truncated fractional maximal function in Lebesgue spaces.
	\begin{lemma}\label{L^pmax}
		Given $R>0$, $1\le p<q<\infty$, $0\le \beta<\frac{n}{q}$ and $f\in L^q(\mathbb{R}^n)$. Then there exists a positive constant $C=C(n,\beta,p,q)>0$ such that
		\begin{equation*}
			\|\mathbf{M}_{\beta,p}[f](\cdot,R)\|_{L^{\frac{nq}{n-\beta q}}(\mathbb{R}^n)}\le C\|f\|_{L^q(\mathbb{R}^n)}.
		\end{equation*}
	\end{lemma}
	\begin{proof}
		For any $x\in \mathbb{R}^n$ and $0<r\le R$, we set the local integral average
		\begin{equation*}
			I_r:=r^\beta\left(\fint_{B_r(x)}|f(y)|^pdy\right)^{1/p}.
		\end{equation*}Then it follows that
		\begin{equation}\label{m1}
			\mathbf{M}_{\beta,p}[f](x,R)=\sup_{0<r\le R}I_r\le \sup_{r>0}I_r.
		\end{equation}Firstly, according to the classical Hardy-Littlewood maximal function, we have
		\begin{equation*}
			I_r\le r^\beta\left(\mathbf{M}[|f|^p](x)\right)^{1/p}.
		\end{equation*}Meanwhile, we apply $q>p$ and H\"{o}lder inequality to yield that
		\begin{equation*}
			\fint_{B_r(x)}|f(y)|^pdy\le |B_1|^{-\frac{p}{q}}r^{-\frac{np}{q}}\|f\|^p_{L^q(\mathbb{R}^n)}.
		\end{equation*}Combining the above estimates, we obtain
		\begin{equation*}
			I_r\le \min\{r^\beta\left(\mathbf{M}[|f|^p](x)\right)^{1/p},C(n,q)r^{\beta-n/q}\|f\|_{L^q(\mathbb{R}^n)}\}.
		\end{equation*} We note that $\min\{a,b\}\le a^{1-\theta}b^\theta$ holds whenever $0\le \theta\le 1$. We further exploit this fact for the above estimate of $I_r$ by choosing the exponent $\theta=\frac{\beta q}{n}\in [0,1)$ to derive
		\begin{equation}\label{m2}
			\begin{aligned}
				I_r&\le C(n,q,\beta)\left(r^\beta\left(\mathbf{M}[|f|^p](x)\right)^{1/p}\right)^{1-\theta}\left(r^{\beta-n/q}\|f\|_{L^q(\mathbb{R}^n)}\right)^\theta\\
				&=C(n,q,\beta)\|f\|_{L^q(\mathbb{R}^n)}^\theta\mathbf{M}[|f|^p](x)^{\frac{1-\theta}{p}}.
			\end{aligned}
		\end{equation}Now inserting \eqref{m2} into \eqref{m1}, we obtain
		\begin{equation*}
				\mathbf{M}_{\beta,p}[f](x,R)\le C(n,q,\beta)\|f\|_{L^q(\mathbb{R}^n)}^\theta\mathbf{M}[|f|^p](x)^{\frac{1-\theta}{p}}.
		\end{equation*}Finally, by virtue of the strong $(q/p,q/p)$ boundedness of the Hardy-Littlewood maximal function, we have
		\begin{equation*}
			\begin{aligned}
				\int_{\mathbb{R}^n}\left(\mathbf{M}_{\beta,p}[f](x,R)\right)^{\frac{nq}{n-\beta q}}dx&\le C\|f\|_{L^q(\mathbb{R}^n)}^{\frac{\beta q^2}{n-\beta q}}\int_{\mathbb{R}^n}\left(\mathbf{M}[|f|^p](x)\right)^{q/p}dx\\
				&\le C\|f\|_{L^q(\mathbb{R}^n)}^{\frac{\beta q^2}{n-\beta q}}\int_{\mathbb{R}^n}\left(|f(x)|^p\right)^{q/p}dx\\
				&\le C\|f\|_{L^q(\mathbb{R}^n)}^{\frac{nq}{n-\beta q}}.
			\end{aligned}
		\end{equation*}Thus, the proof of Lemma \ref{L^pmax} is complete.
	\end{proof}
	Referring to \cite[Proposition 2.8]{DN2025}, a rescaling of the result directly yields the following property, which allows us to obtain estimates in fractional Sobolev spaces $W^{t,q}$ by proving appropriate estimates in terms of fractional sharp maximal functions.
	\begin{lemma}\label{Wtq}
		Let $R>0$, $1<\tilde{q}<q<\infty$ and $\beta,\tilde{\beta}\in (0,1)$. If $\beta-n/q=\tilde{\beta}-n/\tilde{q}$, then for any $g\in L^{\tilde{q}}(B_{2R}(x_0)),$ there exists the constant $C=C(n,\beta,\tilde{\beta},q,\tilde{q})>0$ such that
		\begin{equation*}
			R^{-{\beta}}\|g\|_{L^{{q}}(B_{R/8}(x_0))}+[g]_{W^{\beta,q}(B_{R/8}(x_0))}\le C\left(R^{-\tilde{\beta}}\|g\|_{L^{\tilde{q}}(B_{R/8}(x_0))}+\|\mathbf{M}_{\tilde{\beta}}^{\#}[g](\cdot,R)\|_{L^{\tilde{q}}(B_{R}(x_0))}\right).
		\end{equation*}
	\end{lemma}
	\section{Pointwise and oscillation potential estimates for solutions}\label{section3}
	In this section, we prove pointwise potential estimates (Theorem \ref{thm1.1}) and universal potential estimates (Theorem \ref{thm1.3}) for viscosity solutions to the nonlocal fully nonlinear elliptic equation \eqref{prob1}. For notational convenience, we first define the global excess functional of $u\in C(\Omega)\cap L^1_{2s}(\mathbb{R}^n)$ by
	\begin{equation*}
		E(u;x,r):= A(u;x,r)+\operatorname{Tail}(u-(u)_{B_r(x)};x,r),
	\end{equation*}
	where
\begin{equation*}
	A(u;x,r):= \fint_{B_r(x)}|u(y)-(u)_{B_r(x)}|dy
\end{equation*}
	is the standard local excess functional of $u$, and the ball $B_r(x)\Subset\Omega$.

	As a preliminary step towards the potential estimates for nonhomogeneous equations, we first establish a decay estimate for the global excess functional of the viscosity solution to
	\begin{equation}\label{v}
		\mathcal{I}(v,x)=0\quad\text{in }\Omega,
	\end{equation}where $\mathcal{I}\in \mathfrak{J}^0_s(\lambda,\Lambda)$.
	\begin{lemma}\label{alpha}
		Let $s\in(0,1)$ and $B_R(x_0)\Subset \Omega$ be a ball. Assume that $v\in C(\Omega)\cap L^1_{2s}(\mathbb{R}^n)$ is a viscosity solution to \eqref{v}. Then there exists a positive constant $C=C(n,s,\lambda,\Lambda)$ such that
		\begin{equation*}
			E(v;x_0,\rho)\le C\left(\frac{\rho}{r}\right)^{\alpha_0}E(v;x_0,r)
		\end{equation*}for any $0<\rho\le r$, where $\alpha_0=\alpha_0(n,s,\lambda,\Lambda)$ is the H\"{o}lder continuity exponent from Lemma \ref{holder}.
	\end{lemma}
	\begin{proof}
	By the definition of the nonlocal excess functional of $v$, we estimate the local excess functional and the tail term separately in what follows.
		\step The objective of this step is to prove that the local excess functional of $v$ satisfies the following decay estimate
		\begin{equation}\label{local}
			A(v;x_0,\rho)=\fint_{B_\rho(x_0)}|v(y)-(v)_{B_\rho(x_0)}|dy\le C\left(\frac{\rho}{r}\right)^{\alpha_0}E(v;x_0,r)
		\end{equation}
		for any $0<\rho\le r$, where the positive constant $C$ depends only on $n,s,\lambda$ and $\Lambda$.
		Observe that the above inequality holds trivially when $r/8<\rho\le r$. Without loss of generality, we may assume that $0<\rho\le r/8$. Consider a rescaling function
		\begin{equation*}
			v_r(z):=v(x+rz)
		\end{equation*}for $x\in B_\rho(x_0)$ and $z\in \mathbb{R}^n$. We denote that $\tilde{v_r}(z):=v_r(z)-(v)_{B_r(x_0)}$.
		
		When $z\in B_{1/2}$, we know that $x+rz\in B_{r/2}(x)\subset B_r(x_0)$. Thus by using Lemma \ref{elliptic}, we have
		\begin{equation*}
			\mathcal{P}^{+}\tilde{v_r}(z)=\mathcal{P}^{+}\left(v(x+rz)-(v)_{B_r(x_0)}\right)=r^{2s}\mathcal{P}^{+}v(x+rz)\ge 0\quad\text{in }B_{1/2}
		\end{equation*}in the viscosity sense. We similarly obtain that
		\begin{equation*}
		\mathcal{P}^{-}\tilde{v_r}\le 0\quad\text{in $B_{1/2}$ in the viscosity sense}.
		\end{equation*}
		According to Lemma \ref{holder}, one deduces that $\tilde{v_r}$ has local H\"{o}lder continuity and satisfies
		\begin{equation*}
			\|\tilde{v_r}\|_{C^{\alpha_0}(B_{1/4})}\le C\|\tilde{v_r}\|_{L^1_{2s}(\mathbb{R}^n)}
		\end{equation*}where the positive constant $C=C(n,s,\lambda,\Lambda)$.
		Then it holds that
		\begin{equation*}
			\begin{aligned}
				\fint_{B_\rho(x_0)}|v(y)-(v)_{B_\rho(x_0)}|dy&\le \sup_{B_\rho(x_0)}|v(x)-v(y)|\\
				&\le \sup_{B_\rho(x_0)}\left|\tilde{v_r}(0)-\tilde{v_r}\left(\frac{y-x}{r}\right)\right|\\
				&\le C\left(\frac{\rho}{r}\right)^{\alpha_0}\|\tilde{v_r}\|_{L^1_{2s}(\mathbb{R}^n)}.
			\end{aligned}
		\end{equation*}
		Next, we claim that $\|\tilde{v_r}\|_{L^1_{2s}(\mathbb{R}^n)}\le C(n,s)E(v;x_0,r)$. We divide $\|\tilde{v_r}\|_{L^1_{2s}(\mathbb{R}^n)}$ into two parts as follows
		\begin{equation*}
			\|\tilde{v_r}\|_{L^1_{2s}(\mathbb{R}^n)}=\int_{B_{1/2}}\frac{|\tilde{v_r}(z)|}{1+|z|^{n+2s}}dz+\int_{\mathbb{R}^n\setminus B_{1/2}}\frac{|\tilde{v_r}(z)|}{1+|z|^{n+2s}}dz.
		\end{equation*}
		Firstly, if $z\in B_{1/2}$, then by change of variables we derive
			\begin{align}\label{3_1}
				\int_{B_1}\frac{|\tilde{v_r}(z)|}{1+|z|^{n+2s}}dz&\le \int_{B_1} |v(x+rz)-(v)_{B_r(x_0)}|dz\nonumber\\
				&\le C(n)\fint_{B_r(x_0)}|v(y)-(v)_{B_r(x_0)}|dy.
			\end{align}
		If instead $z\in \mathbb{R}^n\setminus B_{1/2}$, then it holds that
			\begin{align}\label{3_2}
				&\int_{\mathbb{R}^n\setminus B_{1/2}}\frac{|\tilde{v_r}(z)|}{1+|z|^{n+2s}}dz\nonumber\\
				\le& \int_{\mathbb{R}^n\setminus B_{1/2}}\frac{|\tilde{v_r}(z)|}{|z|^{n+2s}}dz\nonumber\\
				=&r^{2s}\int_{\mathbb{R}^n\setminus B_{r/2}(x)}\frac{|v(y)-(v)_{B_r(x_0)}|}{|y-x|^{n+2s}}dy\nonumber\\
				=&r^{2s}\int_{\mathbb{R}^n\setminus B_{r}(x_0)}\frac{|v(y)-(v)_{B_r(x_0)}|}{|y-x|^{n+2s}}dy+r^{2s}\int_{B_r(x_0)\setminus B_{r/2}(x)}\frac{|v(y)-(v)_{B_r(x_0)}|}{|y-x|^{n+2s}}dy\nonumber\\
				\le& C(n,s)\left(r^{2s}\int_{\mathbb{R}^n\setminus B_{r}(x_0)}\frac{|v(y)-(v)_{B_r(x_0)}|}{|y-x_0|^{n+2s}}dy+\fint_{B_r(x_0)}|v(y)-(v)_{B_r(x_0)}|dy\right).
			\end{align}
		Combining \eqref{3_1} with \eqref{3_2}, we finish the proof of claim. Hence, the local decay estimae \eqref{local} is valid for any $0<\rho\le r$.
		
		\step In this step, we turn our attention to establishing the decay estimate for the tail term. For $0<\rho\le r$, we split the tail term into annulus part and far field as follows
		\begin{equation*}
			\begin{aligned}
				&\operatorname{Tail}(v-(v)_{B_\rho(x_0)};x_0,\rho)\\
				=&\rho^{2s}\int_{\mathbb{R}^n\setminus B_{r}(x_0)}\frac{|v(y)-(v)_{B_\rho(x_0)}|}{|y-x_0|^{n+2s}}dy+\rho^{2s}\int_{B_r(x_0)\setminus B_{\rho}(x_0)}\frac{|v(y)-(v)_{B_\rho(x_0)}|}{|y-x_0|^{n+2s}}dy\\
				:=&T_1+T_2.
			\end{aligned}
		\end{equation*}
		For the integral over $\mathbb{R}^n \setminus B_r(x_0)$, we use the triangle inequality to get
		\begin{equation*}
			\begin{aligned}
				T_1&=\left(\frac{\rho}{r}\right)^{2s}\operatorname{Tail}(v-(v)_{B_\rho(x_0)};x_0,r)\\
				&\le \left(\frac{\rho}{r}\right)^{2s}\left(\operatorname{Tail}(v-(v)_{B_r(x_0)};x_0,r)+\operatorname{Tail}((v)_{B_r(x_0)}-(v)_{B_\rho(x_0)};x_0,r)\right)\\
				&\le C\left(\frac{\rho}{r}\right)^{2s}E(v;x_0,r)+C\left(\frac{\rho}{r}\right)^{2s}|(v)_{B_r(x_0)}-(v)_{B_\rho(x_0)}|,
			\end{aligned}
		\end{equation*}where $C=C(n,s)>0$. For estimating the last term of above inequality, we divide the range of $\rho$ into two parts. If $\rho\in (r/4,r]$,
		\begin{equation*}
			\begin{aligned}
				|(v)_{B_r(x_0)}-(v)_{B_\rho(x_0)}|&\le \fint_{B_\rho(x_0)}|v-(v)_{B_r(x_0)}|dx\\
				&\le C(n)\fint_{B_r(x_0)}|v-(v)_{B_r(x_0)}|dx\\
				&\le C(n)E(v;x_0,r) \\
				&\le C(n,s)\left(\frac{\rho}{r}\right)^{2s}E(v;x_0,r).
			\end{aligned}
		\end{equation*}Next we consider the case $\rho\in (0,r/4]$, then there exist $\gamma\in (1/4,1/2]$ and positive integer $k\ge 2$ such that $\rho=\gamma^k r$. By summing up through splitting terms, we have
		\begin{equation*}
			\begin{aligned}
				|(v)_{B_\rho(x_0)}-(v)_{B_{\gamma r}(x_0)}|
				&\le \sum_{j=1}^{k-1}|(v)_{B_{\gamma^{j}r}(x_0)}-(v)_{B_{\gamma^{j+1}r}(x_0)}|\\
				&\le \gamma^{-n}\sum_{j=1}^{k-1}\fint_{B_{\gamma^j r}(x_0)}|v-(v)_{B_{\gamma^{j}r}(x_0)}|dx\\
				&\le C\sum_{j=1}^{k-1}\int_{\gamma^j r}^{\gamma^{j-1}r}A(v;x_0,\gamma^j r)\frac{dt}{t}\\
				&\le C\sum_{j=1}^{k-1}\int_{\gamma^j r}^{\gamma^{j-1}r}A(v;x_0,t)\frac{dt}{t}\\
				&\le C\int_{\rho}^{r}A(v;x_0,t)\frac{dt}{t}.
			\end{aligned}
		\end{equation*}Summarize two cases together, for any $\rho\in (0,r]$, it holds that
		\begin{equation}\label{p2}
			T_1\le C\left(\frac{\rho}{r}\right)^{2s}E(v;x_0,r)+C\int_{\rho}^{r}\left(\frac{\rho}{t}\right)^{2s}A(v;x_0,t)\frac{dt}{t}.
		\end{equation}
		
		For the integral over $B_r(x_0)\setminus B_\rho(x_0)$ and firstly, we consider $\rho\in(0,r/4]$ and obtain that

			\begin{align}\label{p1}
				T_2&=\rho^{2s}\int_{B_r(x_0)\setminus B_{\rho}(x_0)}\frac{|v(y)-(v)_{B_\rho(x_0)}|}{|y-x_0|^{n+2s}}dy\nonumber\\
				&=\rho^{2s}\sum_{j=0}^{k-1}\int_{B_{\gamma^{-j-1}\rho}(x_0)\setminus B_{\gamma^{-j}\rho}(x_0)}\frac{|v(y)-(v)_{B_\rho(x_0)}|}{|y-x_0|^{n+2s}}dy\nonumber\\
				&\le \rho^{2s}\sum_{j=0}^{k-1}\left(\gamma^{-j}\rho\right)^{-n-2s}\left(\gamma^{-j-1}\rho\right)^{n}\fint_{B_{\gamma^{-j-1}\rho}(x_0)}|v-(v)_{B_\rho(x_0)}|dy\nonumber\\
				&\le C(n,s)\sum_{j=0}^{k}\gamma^{2sj}\fint_{B_{\gamma^{-j}\rho}(x_0)}|v-(v)_{B_\rho(x_0)}|dy\nonumber\\
				&\le C(n,s)\sum_{j=0}^{k}\gamma^{2sj}A(v;x_0,\gamma^{-j}\rho)\nonumber\\
				&\le C(n,s)\sum_{j=0}^{k}\int_{\gamma^{-j}\rho}^{\gamma^{-j-1}\rho}\gamma^{2sj}A(v;x_0,\gamma^{-j}\rho)\frac{dt}{t}+C(n,s)\gamma^{2sk}A(v;x_0,r)\nonumber\\
				&\le C(n,s)\int_{\rho}^{r}\left(\frac{\rho}{t}\right)^{2s}A(v;x_0,t)\frac{dt}{t}+C(n,s)\left(\frac{\rho}{r}\right)^{2s}E(v;x_0,r).
			\end{align}
		For $\rho\in (r/4,r]$,
		\begin{equation}\label{p3}
			T_2\le A(v;x_0,r)\le C\left(\frac{\rho}{r}\right)^{2s}E(v;x_0,r)
		\end{equation}
	
		Combining \eqref{p2}, \eqref{p1} with \eqref{p3}, we obtain the desired decay for the tail term
		\begin{equation}
			\operatorname{Tail}(v-(v)_{B_\rho(x_0)};x_0,\rho)\le C\left(\frac{\rho}{r}\right)^{2s}E(v;x_0,r)+C\int_{\rho}^{r}\left(\frac{\rho}{t}\right)^{2s}A(v;x_0,t)\frac{dt}{t}.
		\end{equation}For any $\rho\in (0,r]$, according to \eqref{local}, we have
		\begin{equation*}
			\begin{aligned}
				\int_{\rho}^{r}\left(\frac{\rho}{t}\right)^{2s}A(v;x_0,t)\frac{dt}{t}&\le C\int_{\rho}^{r}\left(\frac{\rho}{t}\right)^{2s}\left(\frac{t}{r}\right)^{\alpha_0}\frac{dt}{t}E(v;x_0,r)\\
				&\le C\left(\frac{\rho}{r}\right)^{\alpha_0}E(v;x_0,r),
			\end{aligned}
		\end{equation*}Since $\alpha_0\le 2s$, finally we obtain that
		\begin{equation*}
			E(v;x_0,\rho)\le C\left(\frac{\rho}{r}\right)^{\alpha_0}E(v;x_0,r).
		\end{equation*}Therefore, we complete the proof of Lemma \ref{alpha}.
	\end{proof}
	To carry over the decay estimate of the global excess functional from the viscosity solution $v$ of the homogeneous equation \eqref{v} to the viscosity solution $u$ of the nonhomogeneous equation \eqref{prob1}, we establish a comparison estimate between $u$ and $v$.
	\begin{lemma}\label{comparison}
		Let $s\in(0,1)$, $p^\ast=p^\ast(n,\lambda,\Lambda)\in [n/2,n)$ and $f\in L^p(\Omega)\cap C(\Omega)$ with $p>\max\{p^{\ast},n/2s\}$. For $B_r(x_0)\Subset \Omega$, assume that $u\in C(\Omega)\cap L^1_{2s}(\mathbb{R}^n)$ is a viscosity solution to \eqref{prob1}, $v\in C(\Omega)\cap L^1_{2s}(\mathbb{R}^n)$ is a viscosity solution to
		\begin{equation}\label{u=v}
			\left\{
			\begin{array}{ll}
				\mathcal{I}(v,x)=0 &\text{in } B_r(x_0),\\
				v=u &\text{in } \mathbb{R}^n\setminus B_r(x_0),
			\end{array}
			\right.
		\end{equation}where $\mathcal{I}\in \mathfrak{J}^0_s(\lambda,\Lambda)$. Then there exists a constant $C=C(n,s,p,\lambda,\Lambda)>0$ such that
		\begin{equation*}
			\|u-v\|_{L^\infty(B_r(x_0))}\le Cr^{2s-\frac{n}{p}}\|f\|_{L^p(B_r(x_0))}.
		\end{equation*}
	\end{lemma}
	\begin{proof}
		Denote that $w(x):=u(x)-v(x)$. By using Lemma \ref{elliptic}, one obtains that
		\begin{equation}\label{p4}
				\mathcal{P}^{+}(w(x))\ge f(x)\ge -|f(x)|\quad\text{and}\quad
				\mathcal{P}^{-}(w(x))\le f(x)\le |f(x)|\quad\text{for }x\in B_r(x_0)
		\end{equation}in the viscosity sense. For $y\in B_1$, set $x:=x_0+ry$, $w_r(y):=w(x_0+ry)$ and $f_r(y):=r^{2s}f(x_0+ry)$. Then
		\begin{equation*}
			\mathcal{P}^{\pm}(w_r(y))=r^{2s}\mathcal{P}^{\pm}(w(x_0+ry))=r^{2s}\mathcal{P}^{\pm}(w(x)).
		\end{equation*}
		According to \eqref{p4}, we have
		\begin{equation*}
				\mathcal{P}^{+}(w_r(y))\ge -|f_r(y)|\quad\text{and}\quad
				\mathcal{P}^{-}(w_r(y))\le |f_r(y)|\quad\text{for }y\in B_1
		\end{equation*}in the viscosity sense and $w_r=0$ in $\mathbb{R}^n\setminus B_1$. By applying Lemma \ref{ABP} to $w_r$, it holds that
		\begin{equation*}
			\|w_r\|_{L^{\infty}(B_1)}\le C\|f_r\|_{L^p(B_1)}.
		\end{equation*}
		Next by rescaling to $w_r$ and $f_r$, we have
		\begin{equation*}
			\|u-v\|_{L^{\infty}(B_r(x_0))}=\|w\|_{L^{\infty}(B_r(x_0))}\le Cr^{2s-\frac{n}{p}}\|f\|_{L^p(B_r(x_0))}.
		\end{equation*}Finally, we finish the proof of Lemma \ref{comparison}.
	\end{proof}

Combining the decay estimate for the global excess functional of the viscosity solution to the homogeneous equation with the preceding comparison estimate, we now deduce the global excess decay for the viscosity solution to the nonhomogeneous equation.
	\begin{lemma}\label{globalexcess}
		Let $s\in(0,1)$, $p^\ast=p^\ast(n,\lambda,\Lambda)\in [n/2,n)$ and $f\in L^p(\Omega)\cap C(\Omega)$ with $p>\max\{p^{\ast},n/2s\}$.
		Assume that $u\in C(\Omega)\cap L^1_{2s}(\mathbb{R}^n)$ is a viscosity solution to \eqref{prob1}, $\mathcal{I}\in \mathfrak{J}^0_s(\lambda,\Lambda)$. Then there exists a positive constant $C=C(n,s,p,\lambda,\Lambda)$ such that
		\begin{equation}\label{decay}
			E(u;x_0,\sigma r)\le C\sigma^{\alpha_0} E(u;x_0,r)+Cr^{2s-\frac{n}{p}}\|f\|_{L^p(B_r(x_0))}
		\end{equation}holds for $B_{\sigma r}(x_0)\subset B_{r}(x_0)\Subset \Omega$, where $\sigma\in(0,1)$ and the exponent $\alpha_0$ is given by Lemma \ref{holder}.
	\end{lemma}
	\begin{proof}
		Let $v$ be the viscosity solution to the homogeneous equation \eqref{u=v} and $w=u-v$. We decompose $u=v+w$ and by triangle inequality, we have
		\begin{equation}\label{sum}
			E(u;x_0,\sigma r)\le E(v;x_0,\sigma r)+E(w;x_0,\sigma r).
		\end{equation}
		We first estimate the the global excess functional $E(v;x_0,\sigma r)$. By using Lemma \ref{alpha} to $v$, we obtain
		\begin{equation}\label{sum1}
			 E(v;x_0,\sigma r)\le C\sigma^{\alpha_0} E(v;x_0,r)\le C\sigma^{\alpha_0}E(u;x_0,r)+C\sigma^{\alpha_0}E(w;x_0,r).
		\end{equation}
		Next, we separately estimate the local excess functional of $w$ and the corresponding tail term as follows
		\begin{equation*}
			\fint_{B_{\sigma r}(x_0)}|w-(w)_{B_{\sigma r}(x_0)}|dx\le 2 \|w\|_{L^{\infty}(B_r(x_0))}
		\end{equation*}and
		\begin{equation*}
			\begin{aligned}
				&(\sigma r)^{2s}\int_{\mathbb{R}^n\setminus B_{\sigma r}(x_0)}\frac{|w(y)-(w)_{B_{\sigma r}(x_0)}|}{|x_0-y|^{n+2s}}dy\\
				\le& (\sigma r)^{2s}\int_{B_r(x_0)\setminus B_{\sigma r}(x_0)}\frac{|w(y)-(w)_{B_{\sigma r}(x_0)}|}{|x_0-y|^{n+2s}}dy+(\sigma r)^{2s}\int_{\mathbb{R}^n\setminus B_{r}(x_0)}\frac{\|w\|_{L^\infty(B_r(x_0))}}{|x_0-y|^{n+2s}}dy\\
				\le& C\|w\|_{L^\infty(B_r(x_0))}.
			\end{aligned}
		\end{equation*}Combining Lemma \ref{comparison} with the above estimates, we derive
		\begin{equation*}
			E(w;x_0,\sigma r)\le C r^{2s-\frac{n}{p}}\|f\|_{L^p(B_r(x_0))}.
		\end{equation*}By a similar argument, we also have
		\begin{equation*}
			E(w;x_0,r)\le Cr^{2s-\frac{n}{p}}\|f\|_{L^p(B_r(x_0))}.
		\end{equation*}
		Substituting into \eqref{sum} and \eqref{sum1}, the proof of Lemma \ref{globalexcess} is immediately complete.
	\end{proof}
	We are now in a position to prove our first main result (i.e., Theorem \ref{thm1.1}) concerning the pointwise potential estimate for the viscosity solution of the nonlocal fully nonlinear elliptic equation \eqref{prob1}.
	\begin{proof}[Proof of Theorem \ref{thm1.1}]
		Given $B_R(x)\Subset\Omega$ and set $0<r\le R/4$. According to Lemma \ref{globalexcess}, we know that for any $0<\rho\le r$
		\begin{equation*}
			E(u;x,\sigma\rho)\le C\sigma^{\alpha_0} E(u;x,\rho)+C\rho^{2s-\frac{n}{p}}\|f\|_{L^p(B_\rho(x))}.
		\end{equation*}By choosing $\sigma=\sigma(n,s,p,\lambda,\Lambda) \in (0, 1/4]$ small enough such that $C \sigma^{\alpha_0} \le 1/2$, we sum the iterative excess bounds over all scales as follows.
		Firstly, we consider a sequence of small ball with the same center $x$ and theirs radius are defined by $r_j:=\sigma^jr$. Next we have
		\begin{align}\label{1.1.1}
				|u(x)|&\le|(u)_{B_r(x)}|+\sum_{j=0}^{\infty}\left|(u)_{B_{r_{j+1}}(x)}-(u)_{B_{r_{j}}(x)}\right|\nonumber\\
				&\le |(u)_{B_r(x)}|+\sigma^{-n}\sum_{j=0}^{\infty}\fint_{B_{r_{j}}(x)}|u(z)-(u)_{B_{r_{j}}(x)}|dz\nonumber\\
				&\le C\fint_{B_r(x)}|u(y)|dy+\sigma^{-n}\sum_{j=0}^{\infty}E(u;x,r_j).
			\end{align}By applying Lemma \ref{globalexcess}, we have
		\begin{equation*}
			E(u;x,r_{j+1})\le \frac{1}{2}E(u;x,r_j)+Cr_{j}^{2s-\frac{n}{p}}\|f\|_{L^p(B_{r_j}(x))} \quad\text{for }j\in\mathbb{N}.
		\end{equation*}Then by summing from $j=0$ to $N$, we obtain
		\begin{equation*}
			\begin{aligned}
				\sum_{j=0}^{N}E(u;x,r_{j+1})&\le \frac{1}{2}\sum_{j=0}^{N}E(u;x,r_j)+C\sum_{j=0}^{N}r_{j}^{2s-\frac{n}{p}}\|f\|_{L^p(B_{r_j}(x))},
			\end{aligned}
		\end{equation*}
		Then it follows that
		\begin{equation*}
			\sum_{j=0}^{N+1}E(u;x,r_{j})\le 2E(u;x,r_0)+C\sum_{j=0}^{N}r_{j}^{2s-\frac{n}{p}}\|f\|_{L^p(B_{r_j}(x))}.
		\end{equation*}Let $N\to \infty$, we further derive
		\begin{equation}\label{1.1.2}
			\sum_{j=0}^{\infty}E(u;x,r_{j})\le CE(u;x,r)+C\sum_{j=0}^{\infty}r_{j}^{2s-\frac{n}{p}}\|f\|_{L^p(B_{r_j}(x))}.
		\end{equation}Since
		\begin{equation*}
			\|f\|_{L^p(B_{r_j}(x))}r_j^{2s-\frac{n}{p}}\le C\int_{r_{i}}^{r_{i-1}}\frac{t^{2s}\left(\fint_{B_t(x)}|f(y)|^pdy\right)^{1/p}}{t}dt,
		\end{equation*}one obtains
		\begin{equation}\label{1.1.3}
			\sum_{j=0}^{\infty}r_{j}^{2s-\frac{n}{p}}\|f\|_{L^p(B_{r_j}(x))}\le C\mathbf{I}^{f}_{2s,p}(x,\frac{r}{\sigma}).
		\end{equation}According to \eqref{1.1.2} and \eqref{1.1.3}, we set $r= \sigma R\le R/4$ and obtain that
\begin{equation*}
	\sum_{j=0}^{\infty}E(u;x,r_{j})\le  C\left(E(u;x,R)+\mathbf{I}^{f}_{2s,p}(x,R)\right).
\end{equation*}Finally, combining the above equation with \eqref{1.1.1}, we deduce that
		\begin{equation*}
			|u(x)|\le C\left(\mathbf{I}^{f}_{2s,p}(x,R)+\fint_{B_R(x)}|u(y)|dy+\operatorname{Tail}(u;x,R)\right).
		\end{equation*}Therefore, we complete the proof of Theorem \ref{thm1.1}.
	\end{proof}

	To further prove the universal potential estimates stated in Theorem \ref{thm1.3}, we first establish a pointwise estimate for the truncated nonlocal fractional sharp maximal function of the viscosity solution to equation \eqref{prob1}, defined by
	\begin{equation*}
		\mathbf{N}_{t}^{\#}[u](x_0,R):=\sup_{0<r\le R}r^{-t}E(u;x_0,r),
	\end{equation*}where the exponent $0\le t<1$ and the ball $B_R(x_0)\Subset\Omega$.
	\begin{theorem}\label{maximal}
		Let $s\in(0,1)$, $p^\ast=p^\ast(n,\lambda,\Lambda)\in [n/2,n)$ and $f\in L^p(\Omega)\cap C(\Omega)$ with $p>\max\{p^{\ast},n/2s\}$. Assume that $u\in C(\Omega)\cap L^1_{2s}(\mathbb{R}^n)$ is a viscosity solution to \eqref{prob1} of the form \eqref{i}. For any $\tilde{\alpha}\in (0,\alpha_0)$, there exists a positive constant $C=C(n,s,p,\lambda,\Lambda,\tilde{\alpha})$ such that
		\begin{equation*}
			\mathbf{N}_{t}^{\#}[u](x_0,R)\le C\left(\mathbf{M}_{2s-t,p}[f](x_0,R)+R^{-t}E(u;x_0,R)\right),
		\end{equation*}for any $0\le t\le\tilde{\alpha}$ and any $B_R(x_0)\Subset\Omega$,
		where  $\alpha_0=\alpha_0(n,s,\lambda,\Lambda)$ is the H\"{o}lder continuity exponent given by Lemma \ref{holder}.
	\end{theorem}
	\begin{proof}
		By applying \eqref{decay} in Lemma \ref{globalexcess} and multiplying both sides by $(\sigma\rho)^{-t}$ with $0<\rho\le R$, we know
		\begin{equation*}
			(\sigma\rho)^{-t}E(u;x_0,\sigma\rho)\le C\sigma^{\alpha_0-t}\rho^{-t}E(u;x_0,\rho)+C\sigma^{-t}\rho^{2s-\frac{n}{p}-t}\|f\|_{L^p(B_\rho)}.
		\end{equation*} We can choose $\sigma=\sigma(n,s,p,\lambda,\Lambda,\tilde{\alpha})>0$ small enough satisfying that $C\sigma^{\alpha_0-t}\le C\sigma^{\alpha_0-\tilde{\alpha}} \le 1/2$, then
		\begin{equation*}
			\begin{aligned}
				(\sigma\rho)^{-t}E(u;x_0,\sigma\rho)&\le \frac{1}{2}\rho^{-t}E(u;x_0,\rho)+C\rho^{2s-\frac{n}{p}-t}\|f\|_{L^p(B_\rho)}\\
				&\le \frac{1}{2}\mathbf{N}_{t}^{\#}[u](x_0,R)+C\mathbf{M}_{2s-t,p}[f](x_0,R).
			\end{aligned}
		\end{equation*}
		Since $0<\rho\le R$ is arbitrarily chosen, we obtain
		\begin{equation}\label{n1}
			\sup_{0<r\le \sigma R}r^{-t}E(u;x_0,r)\le \frac{1}{2}\mathbf{N}_{t}^{\#}[u](x_0,R)+C\mathbf{M}_{2s-t,p}[f](x_0,R).
		\end{equation}On the other hand, by direct integral calculating, we have
		\begin{equation}\label{n2}
			\sup_{\sigma R<r\le R}r^{-t}E(u;x_0,r)\le CR^{-t}E(u;x_0,R).
		\end{equation}
		Combining \eqref{n1} with \eqref{n2} and absorbing the first term on the right-hand side of \eqref{n1} into the left-hand side, then we complete the proof of Theorem \ref{maximal}.	
	\end{proof}
This section concludes with a proof of the oscillation potential  estimate for viscosity solutions to the nonlocal fully nonlinear elliptic equation \eqref{prob1}.

	\begin{proof}[Proof of Theorem \ref{thm1.3}]
		Given $x\in B_{R/8}(x_0)$, we set $r\le R/4$ to ensure $B_r(x)\subset B_R(x_0)$. Basis on Lemma \ref{globalexcess}, we know that for any $0<\rho\le r$,
		\begin{equation*}
			E(u;x,\sigma\rho)\le C\sigma^{\alpha_0} E(u;x,\rho)+C\rho^{2s-\frac{n}{p}}\|f\|_{L^p(B_\rho(x))}.
		\end{equation*}
		For every $\tilde{\alpha}\in(0,\alpha_0)$, we choose $\sigma=\sigma(n,s,p,\lambda,\Lambda,\tilde{\alpha})>0$ small enough such that $C\sigma^{\alpha_0}\le \sigma^{\tilde{\alpha}}\le \sigma^{\alpha}$ with any $\alpha\in[0,\tilde{\alpha}]$, then
		\begin{equation}\label{iter}
			E(u;x,\sigma\rho)\le\sigma^{\alpha} E(u;x,\rho)+C\rho^{2s-\frac{n}{p}}\|f\|_{L^p(B_\rho(x))}.
		\end{equation}

		Next we consider the Campanato iteration in a sequence of small balls with same center. Given initial radius $r>0$, we define the radius of ball sequence by $r_i=\sigma^{i}r$, $i\in\mathbb{N}$. According to \eqref{iter}, we deduce for any integer $k\ge 1$,
		\begin{equation}\label{sequence}
				E(u;x,r_k)\le\sigma^{\alpha k} E(u;x,r)+C\sum_{i=0}^{k-1}\sigma^{\alpha(k-1-i)}r_i^{2s-\frac{n}{p}}\|f\|_{L^p(B_{r_i}(x))}.
		\end{equation}
		Since $\sigma^{k-1-i}=\frac{\sigma^{k-1}}{\sigma^{i}}=\frac{r_{k-1}}{r_{i}}$, we obtain
		\begin{equation*}
			E(u;x,r_k)\le \left(\frac{r_k}{r}\right)^{\alpha} E(u;x,r)+Cr_{k-1}^{\alpha}\sum_{i=0}^{k-1}r_i^{2s-\alpha-\frac{n}{p}}\|f\|_{L^p(B_{r_i}(x))}.
		\end{equation*}
		For $t\in[r_{i},r_{i-1}]$, since $t^{2s}\|f\|_{L^p(B_{t}(x))}$ is increasing about $t$. Then
		\begin{equation*}
			\begin{aligned}
				r_i^{2s-\alpha-\frac{n}{p}}\|f\|_{L^p(B_{r_i}(x))}&=\frac{\alpha+\frac np}{1-\sigma^{\alpha+\frac np}} r_i^{2s}\|f\|_{L^p(B_{r_i}(x))}\int_{r_{i}}^{r_{i-1}}t^{-\alpha-1-\frac{n}{p}}dt\\
				&\le C\int_{r_{i}}^{r_{i-1}}t^{2s-\alpha-\frac{n}{p}-1}\|f\|_{L^p(B_{t}(x))}dt\\
				&\le 	C\int_{r_{i}}^{r_{i-1}}t^{2s-\alpha-1}\left(\fint_{B_t(x)}|f(y)|^pdy\right)^{1/p}dt.
			\end{aligned}
		\end{equation*}
		Summing over $i$ from $0$ to $k-1$, we derive
		\begin{equation*}
			\sum_{i=0}^{k-1}r_i^{2s-\alpha-\frac{n}{p}}\|f\|_{L^p(B_{r_i}(x))}\le C\mathbf{I}^{f}_{2s-\alpha,p}(x,\frac{r}{\sigma}).
		\end{equation*}
		Then for any $0<\rho\le r$, there holds that
		\begin{equation*}
			E(u;x,\rho)\le C\left(\frac{\rho}{r}\right)^{\alpha} E(u;x,r)+C\rho^{\alpha}\mathbf{I}^{f}_{2s-\alpha,p}(x,\frac{r}{\sigma}).
		\end{equation*}We set $r= \sigma R\le R/4$ and consider $0<\rho\le \sigma R$,
		\begin{equation*}
			\begin{aligned}
				E(u;x,\rho)&\le C\left(\frac{\rho}{\sigma R}\right)^{\alpha} E(u;x,\sigma R)+C\rho^{\alpha}\mathbf{I}^{f}_{2s-\alpha,p}(x,R)\\
				&\le C\left(\frac{\rho}{R}\right)^{\alpha} E(u;x, R)+C\rho^{\alpha}\mathbf{I}^{f}_{2s-\alpha,p}(x,R)
			\end{aligned}
		\end{equation*}and for $\sigma R<\rho\le R$, it is easy to verify that
	\begin{equation*}
		E(u;x,\rho)\le C\left(\frac{\rho}{R}\right)^{\alpha}E(u;x,R).
	\end{equation*}Overall, we conclude that
	\begin{equation}\label{iteration}
		E(u;x,\rho)\le C\left(\frac{\rho}{R}\right)^{\alpha}E(u;x,R)+C\rho^{\alpha}\mathbf{I}^{f}_{2s-\alpha,p}(x,R)
	\end{equation}for any $0<\rho\le R$.
Now given $x,y\in B_{R/8}(x_0)$. Note $d:=|x-y|$ and $d_j:=\frac{d}{2^j}\le r$, $j\in\mathbb{N}$. According to Lebesgue differential theorem, given that $u$ is continuous, we conclude that
		\begin{equation}\label{lim}
			u(x)=\lim\limits_{j\to \infty}(u)_{B_{d_j}(x)}\quad\text{for every }x\in B_{R/8}(x_0).
		\end{equation}
		Then it follows that
		\begin{equation*}
			\begin{aligned}
				|u(x)-(u)_{B_d(x)}|&\le \sum_{j=0}^{\infty}|(u)_{B_{d_{j+1}}(x)}-(u)_{B_{d_{j}}(x)}|\\
				&\le \sum_{j=0}^{\infty}\fint_{B_{d_{j+1}}(x)}|u(z)-(u)_{B_{d_{j}}(x)}|dz\\
				&\le \sum_{j=0}^{\infty}2^n\fint_{B_{d_{j}}(x)}|u(z)-(u)_{B_{d_{j}}(x)}|dz\\
				&\le C\sum_{j=0}^{\infty}E(u;x,d_j)\quad\quad\text{for every }x\in B_{R/8}(x_0).
			\end{aligned}
		\end{equation*}We apply \eqref{iteration} to $E(u;x,d_j)$ for each $0<d_j\le r$, and since $\sum\limits_{j=0}^{\infty}(d_j)^\alpha$ is convergence, we obtain
		\begin{equation}\label{part1}
			\begin{aligned}
				|u(x)-(u)_{B_d(x)}|&\le Cd^{\alpha}\left(r^{-\alpha}E(u;x,r)+\mathbf{I}^{f}_{2s-\alpha,p}(x,r)\right).
			\end{aligned}
		\end{equation}
	Next, we introduce a ball $B_{2d}(x)$ that can contain both $B_d(x)$ and $B_d(y)$ and insert it into the inequality
			\begin{align}\label{part2}
				|(u)_{B_d(x)}-(u)_{B_d(y)}|&\le |(u)_{B_d(x)}-(u)_{B_{2d}(x)}|+|(u)_{B_{2d}(x)}-(u)_{B_d(y)}|\nonumber\\
				&\le \fint_{B_d(x)}|u(z)-(u)_{B_{2d}(x)}|dz+\fint_{B_d(y)}|u(z)-(u)_{B_{2d(x)}}|dz\nonumber\\
				&\le CE(u;x,2d)\nonumber\\
				&\le C\left(\frac{d}{r}\right)^{\alpha} E(u;x,r)+Cd^{\alpha}\mathbf{I}^{f}_{2s-\alpha,p}(x,r).
			\end{align}By the triangle inequality, we have
		\begin{equation*}
			|u(x)-u(y)|\le |u(x)-(u)_{B_d(x)}|+|(u)_{B_d(x)}-(u)_{B_d(y)}|+|u(y)-(u)_{B_d(y)}|.
		\end{equation*}Since \eqref{part1} and \eqref{part2}, ones obtain that for any $x,y\in B_{R/8}(x_0)$,
		\begin{equation*}
			|u(x)-u(y)|\le Cd^\alpha\left(\mathbf{I}^{f}_{2s-\alpha,p}(x,r)+\mathbf{I}^{f}_{2s-\alpha,p}(y,r)\right)+C\left(\frac{d}{r}\right)^{\alpha} \left(E(u;x,r)+E(u;y,r)\right).
		\end{equation*}
		Since $B_r(x),B_r(y)\subset B_{2r}(x)$ and Theorem \ref{maximal}, for $r\le R/4$, we have
		\begin{equation*}
			\begin{aligned}
				E(u;x,r)+E(u;y,r)&\le CE(u;x,2r)\\
				&\le CR^\alpha\mathbf{N}_{\alpha}^{\#}[u](x,R/2)\\
				&\le CR^\alpha\left(\mathbf{M}_{2s-\alpha,p}[f](x,R/2)+R^{-\alpha}E(u;x,R/2)\right).
			\end{aligned}
		\end{equation*}By Lemma \ref{Riesz}, we have
		\begin{equation*}
			\mathbf{M}_{2s-\alpha,p}[f](x,R/2)\le C\mathbf{I}^f_{2s-\alpha,p}(x,R).
		\end{equation*}Moreover, since $B_{R/2}(x)\subset B_R(x_0)$, one readily obtains by direct calculation that
		\begin{equation*}
			E(u;x,R/2)\le CE(u;x_0,R)\le C \left( \fint_{B_R(x_0)} |u| dx + \text{Tail}(u; x_0, R) \right).
		\end{equation*}
		Hence, we complete the proof of Theorem \ref{thm1.3}.
	\end{proof}
	
	\section{Gradient potential estimates and higher Sobolev regularity}\label{section4}
	In this section, we develop the potential theory for gradients of viscosity solutions to the nonlocal fully nonlinear elliptic equation \eqref{prob1}. Note that such viscosity solutions may not possess weak derivatives, and the comparison estimate with the viscosity solution  to the corresponding homogeneous equation holds only in the $L^\infty$ norm. For this reason, we introduce a novel class of global excess functionals. To begin with, we define the space consisting of all affine functions on the ball.
	\begin{definition}
		Let $s\in (1/2,1)$ and $B_r(x_0)\Subset \Omega$. The affine functions space in $B_r(x_0)$ is defined as
		\begin{equation*}
			\mathbf{L}_{x_0,r}(\mathbb{R}^n):=\{l(x)=a+A\cdot(x-x_0):(a,A)\in \mathbb{R}\times\mathbb{R}^n,x\in B_r(x_0)\},
		\end{equation*}and the norm
		\begin{equation*}
			\|l\|_{\mathbf{L}_{x_0,r}(\mathbb{R}^n)}:=\frac{|a|}{r}+|A|.
		\end{equation*}
	\end{definition}
It is easy to check that $\mathbf{L}_{x_0,r}(\mathbb{R}^n)\subset L^1_{2s}(\mathbb{R}^n)$ whenever $s>1/2$. We now define the global excess functional in terms of affine function.
	\begin{definition}
		Let $s\in (1/2,1)$ and $B_r(x_0)\Subset\Omega$.
	For any $u\in C(\Omega)\cap L^1_{2s}(\mathbb{R}^n)$, the global excess functional of $u$ associated with affine functions is defined as
		\begin{equation*}
			G(u;x_0,r):=\inf_{l\in \mathbf{L}_{x_0,r}(\mathbb{R}^n)}\mathcal{E}_{l}(u;x_0,r),
		\end{equation*}where
		\begin{equation*}
		\mathcal{E}_{l}(u;x_0,r):=\frac{1}{r}\left(\fint_{B_{r}(x_0)}|u(y)-l(y)|dy+\operatorname{Tail}(u-l;x_0,r)\right).
	\end{equation*}
	\end{definition}
	We need to verify the attainability of infimum in the definition of $G(u;x_0,r)$. For this purpose, we first prove the compatibility of the norms.
	\begin{lemma}\label{compate}
		Let $s\in (1/2,1)$ and $B_r(x_0)\Subset\Omega$. For every $l\in \mathbf{L}_{x_0,r}(\mathbb{R}^n)$ and note $l(x)=a+A\cdot(x-x_0)$, there exist constants $c=c(n,s)>0$ and $C=C(n,s)>0$ such that
		\begin{equation*}
			c\left(\frac{|a|}{r}+|A|\right)\le \mathcal{E}_{l}(0;x_0,r)\le C\left(\frac{|a|}{r}+|A|\right).
		\end{equation*}
	\end{lemma}
\begin{proof}
	 Through variable substitution, we have
	 \begin{equation*}
	 	\fint_{B_{r}(x_0)}|l(y)|dy=\fint_{B_1}|a+rA\cdot z|dz.
	 \end{equation*}Due to the symmetry of ball $B_1$ and $2|a|\le |a+rA\cdot z|+|a-rA\cdot z|$,
 	\begin{equation}\label{eq3}
 		\fint_{B_1}|a+rA\cdot z|dz\ge |a|
 	\end{equation}On the other hand, by $2r|A\cdot z|\le |a+rA\cdot z|+|a-rA\cdot z|$,
 	\begin{equation}\label{eq4}
 		\fint_{B_1}|a+rA\cdot z|dz\ge r\fint_{B_1}|A\cdot z|dz=c(n)r|A|.
 	\end{equation}From \eqref{eq3} and \eqref{eq4}, we obtain the lower bound
 	\begin{equation*}
 		\frac{1}{r}\fint_{B_{r}(x_0)}|l(y)|dy\ge c(n)\left(\frac{|a|}{r}+|A|\right).
 	\end{equation*}

 	Next we consider the upper bound of $\mathcal{E}_{0;x_0,r}(l)$. For local term, by direct calculation we have
 	\begin{equation*}
 		\frac{1}{r}\fint_{B_r(x_0)}|l(y)|dy\le \frac{|a|}{r}+|A| ,
 	\end{equation*}and for tail term, set $z=y-x_0$,
 	\begin{equation*}
 		\frac{1}{r}\operatorname{Tail}(l;x_0,r)\le r^{2s-1}\int_{|z|\ge r}\frac{|a|+|A||z|}{|z|^{n+2s}}dz\le C(n,s)\left(\frac{|a|}{r}+|A|\right).
 	\end{equation*}In summary, it can be concluded that the norms are compatible.
\end{proof}
We are now ready to show that the infimum is attained.
\begin{lemma}\label{reach}
	Let $s\in (1/2,1)$ and $u\in C(\Omega)\cap L^1_{2s}(\mathbb{R}^n)$. For any $x_0\in \mathbb{R}^n$ and $r>0$, there exists an affine function $l_{x_0,r}\in \mathbf{L}_{x_0,r}(\mathbb{R}^n)$ such that
	\begin{equation*}
		G(u;x_0,r)=\mathcal{E}_{l_{x_0,r}}(u;x_0,r).
	\end{equation*}
\end{lemma}
\begin{proof}
	Since $u\in C(\Omega)\cap L^1_{2s}(\mathbb{R}^n)$, we have
	\begin{equation*}
		0\le G(u;x_0,r)\le \mathcal{E}_{0}(u;x_0,r)<\infty.
	\end{equation*}For any $l_1,l_2\in \mathbf{L}_{x_0,r}(\mathbb{R}^n)$, according to Lemma \ref{compate} and $\left||u-l_1|-|u-l_2|\right|\le |l_1-l_2|$, we can obtain
	\begin{equation}
		|\mathcal{E}_{l_1}(u;x_0,r)-\mathcal{E}_{l_2}(u;x_0,r)|\le \mathcal{E}_{l_1-l_2}(0;x_0,r)\le C(n,s)\|l_1-l_2\|_{\mathbf{L}_{x_0,r}(\mathbb{R}^n)}.
	\end{equation}Therefore, $\mathcal{E}_{l}(u;x_0,r)$ is continuous with respect to $l$.

	In addition, from the triangle inequality, we have
	\begin{equation*}
		\mathcal{E}_{l}(u;x_0,r)\ge \mathcal{E}_{l}(0;x_0,r)-\mathcal{E}_{0}(u;x_0,r),
	\end{equation*}again by Lemma \ref{compate}, we further obtain
	\begin{equation*}
		\mathcal{E}_{l}(u;x_0,r)\ge c(n)\left(\frac{|a|}{r}+|A|\right)-\mathcal{E}_{0}(u;x_0,r).
	\end{equation*}As $\frac{|a|}{r}+|A|\to \infty$, it follows that $\mathcal{E}_{l}(u;x_0,r)\to \infty$. That is, $\mathcal{E}_{l}(u;x_0,r)$ satisfies coercivity.

	Since $\mathbf{L}_{x_0,r}(\mathbb{R}^n)\simeq\mathbb R^{n+1}$ is a reflexive Banach space and the functional $\mathcal{E}_l(u;x_0,r)$ is continuous and coercive on $\mathbf{L}_{x_0,r}(\mathbb{R}^n)$, the coercive version of the Weierstrass theorem implies that $\mathcal{E}_l(u;x_0,r)$ attains its minimum. Hence there exists at least one $l_{x_0,r}\in \mathbf{L}_{x_0,r}(\mathbb{R}^n)$ such that
	\begin{equation*}
		G(u;x_0,r)=\mathcal{E}_{l_{x_0,r}}(u;x_0,r).
	\end{equation*}Finally, the proof of Lemma \ref{reach} is complete.
\end{proof}
	Next, we prove that the global excess functional of the new construction satisfies the following Campanato-type decay estimate.
	\begin{lemma}\label{Duiter}Let $s\in(1/2,1)$,  $\sigma\in(0,1)$, $p^\ast=p^\ast(n,\lambda,\Lambda)\in [n/2,n)$ and $f\in L^p(\Omega)\cap C(\Omega)$ with $p>\max\{p^{\ast},n/2s\}$. Assume that $u\in C(\Omega)\cap L^1_{2s}(\mathbb{R}^n)$ is a viscosity solution to $\mathcal{I}(u,x)=f$ in $\Omega$ of the form \eqref{i} with $\mathcal{I}\in \mathfrak{J}^1_s(\lambda,\Lambda)$. Then there exists a positive constant $C=C(n,s,p,\lambda,\Lambda,\omega(\cdot))>0$ such that
			\begin{equation*}
				G(u;x_0,\sigma r)\le C\left(\sigma^{\bar{\alpha}}+\sigma^{2s-1}\right) G(u;x_0,r)+C\sigma^{-1} r^{2s-\frac{n}{p}-1}\|f\|_{L^p(B_r(x_0))},
			\end{equation*}for $B_{\sigma r}(x_0)\subset B_{r}(x_0)\Subset \Omega$, where	$\bar{\alpha}=\bar{\alpha}(n,s,\lambda,\Lambda)\in (0,\alpha_0)$ is the exponent given by Lemma \ref{1+alpha}.
	\end{lemma}
\begin{proof}
	According to Lemma \ref{reach}, there exists an affine function $l_{x_0,r}(x)=a_{x_0,r}+A_{x_0,r}\cdot (x-x_0)$ such that
	\begin{equation*}
		G(u;x_0,r)=\mathcal{E}_{l_{x_0,r}}(u;x_0,r).
	\end{equation*}Since $l_{x_0,r}(x+y)+l_{x_0,r}(x-y)-2l_{x_0,r}(x)=0$, it is easy to know that $\mathcal{I}(u)=\mathcal{I}(u-l_{x_0,r})=f$ in $B_r(x_0)$ in viscosity sense and
	\begin{equation*}
		G(u;x_0,r)=\frac{1}{r}\left(\fint_{B_r(x_0)}|u-l_{x_0,r}|dy+\operatorname{Tail}(u-l_{x_0,r};x_0,r)\right).
	\end{equation*}
	Next we consider the homogeneous equation
	\begin{equation*}
		\left\{
		\begin{array}{ll}
			\mathcal{I}(h,x)=0 & \text{in } B_r(x_0), \\
			h=u-l_{x_0,r} &\text{in } \mathbb{R}^n\setminus B_r(x_0),
		\end{array}
		\right.
	\end{equation*}
	Obviously,  $\operatorname{Tail}(h;x_0,r)=\operatorname{Tail}(u-l_{x_0,r};x_0,r)$ and by applying Lemma \ref{comparison}, we obtain
	\begin{equation}\label{v-h}
		\|u-l_{x_0,r}-h\|_{L^{\infty}(B_r(x_0))}\le Cr^{2s-\frac{n}{p}}\|f\|_{L^p(B_r(x_0))}.
	\end{equation}
	 According to Lemma \ref{1+alpha}, for some $\bar{\alpha}\in (0,\alpha_0)$, $\alpha_0=\alpha_0(n,s,\lambda,\Lambda)$ is the H\"{o}lder continuity exponent in Lemma \ref{holder}, we know that $h\in C^{1,\bar{\alpha}}_{loc}(B_r(x_0))$ and denote $h_r(x):=h(rx+x_0)$,
	\begin{equation}\label{h_1}
		\|h\|_{C^{1,\bar{\alpha}}(B_{r/2}(x_0))}\le Cr^{-1-\bar{\alpha}}\|h_r\|_{L^1_{2s}(\mathbb{R}^n)}.
	\end{equation}
	Consider the affine function $\ell_h(x)=\nabla h(x_0)\cdot (x-x_0)+h(x_0)$, for any $\rho\le r/2$ we have
	\begin{equation*}
		\begin{aligned}
			\fint_{B_\rho(x_0)}|h-\ell_h|dy&\le \|h-\ell_h\|_{L^\infty(B_\rho(x_0))}\\
			&\le C\rho^{1+\bar{\alpha}}r^{-1-\bar{\alpha}}\|h_r\|_{L^1_{2s}(\mathbb{R}^n)}\\
			&\le C\rho^{1+\bar{\alpha}}r^{-1-\bar{\alpha}}\left(\fint_{B_r(x_0)}|h|dy+\operatorname{Tail}(h;x_0,r)\right)\\
			&\le C\rho^{1+\bar{\alpha}}r^{-\bar{\alpha}}G(u;x_0,r)+C\rho^{1+\bar{\alpha}}r^{-1-\bar{\alpha}}\left(r^{2s-\frac{n}{p}}\|f\|_{L^p(B_r(x_0))}\right).
		\end{aligned}
	\end{equation*}
	Let $\rho=\sigma r$ where $\sigma\in(0,1/2)$ is a constant to be determined. Then
		\begin{align}\label{h_2}
			\frac{1}{\sigma r}\fint_{B_{\sigma r}(x_0)}|u-l_{x_0,r}-\ell_h|dy&\le \frac{1}{\sigma r}\fint_{B_{\sigma r}(x_0)}|u-l_{x_0,r}-h|dy+\frac{1}{\sigma r}\fint_{B_{\sigma r}(x_0)}|h-\ell_h|dy\nonumber\\
			&\le C\left(\sigma^{-1}+\sigma^{\bar{\alpha}}\right)r^{2s-\frac{n}{p}-1}\|f\|_{L^p(B_r(x_0))}+C\sigma^{\bar{\alpha}}G(u;x_0,r).
		\end{align}Next we deal with the tail term,
		\begin{align*}
			&\frac{1}{\sigma r}\operatorname{Tail}(u-l_{x_0,r}-\ell_h;x_0,\sigma r)\\
			=&(\sigma r)^{2s-1}\int_{B_{r}(x_0)\setminus B_{\sigma r}(x_0)}\frac{|u-l_{x_0,r}-\ell_h|}{|x_0-y|^{n+2s}}dy+(\sigma r)^{2s-1}\int_{\mathbb{R}^n\setminus B_{r}(x_0)}\frac{|u-l_{x_0,r}-\ell_h|}{|x_0-y|^{n+2s}}dy.
		\end{align*}
	For $|y-x_0|\ge r$ and $\ell_h(y)=\nabla h(x_0)\cdot (y-x_0)+h(x_0)$, there exists a positive universal constant $C=C(n,s)$ such that
		\begin{align}\label{h}
			&\int_{\mathbb{R}^n\setminus B_{r}(x_0)}\frac{|\ell_h(y)|}{|x_0-y|^{n+2s}}dy\nonumber\\
			\le& |\nabla h(x_0)|\int_{\mathbb{R}^n\setminus B_{r}(x_0)}\frac{1}{|x_0-y|^{n+2s-1}}dy+|h(x_0)|\int_{\mathbb{R}^n\setminus B_{r}(x_0)}\frac{1}{|x_0-y|^{n+2s}}dy\nonumber\\
			\le& C\left(r^{-2s+1}|\nabla h(x_0)|+r^{-2s}|h(x_0)|\right),
		\end{align}by using $2s>1$.
	Then according to \eqref{v-h} \eqref{h_1} and \eqref{h}, we have
		\begin{align*}
			&(\sigma r)^{2s-1}\int_{\mathbb{R}^n\setminus B_{r}(x_0)}\frac{|u-l_{x_0,r}-\ell_h|}{|x_0-y|^{n+2s}}dy\\
			\le& \sigma^{2s-1}r^{-1}\operatorname{Tail}(u-l_{x_0,r};x_0,r)+ C\left(\sigma^{2s-1}|\nabla h(x_0)|+\sigma^{2s-1}r^{-1}|h(x_0)|\right)\\
			\le& \sigma^{2s-1}r^{-1}\operatorname{Tail}(u-l_{x_0,r};x_0,r)+C\sigma^{2s-1}r^{-1}\left(\fint_{B_{r}(x_0)}|h|dy+\operatorname{Tail}(h;x_0,r)\right)\\
			\le& C\sigma^{2s-1}r^{-1}\operatorname{Tail}(u-l_{x_0,r};x_0,r)+C\sigma^{2s-1}r^{-1}\left(\fint_{B_{r}(x_0)}|u-l_{x_0,r}|dy+\fint_{B_{r}(x_0)}|u-l_{x_0,r}-h|dy\right)\\
			\le& C\sigma^{2s-1}G(u;x_0,r)+C\sigma^{2s-1}r^{2s-\frac{n}{p}-1}\|f\|_{L^p(B_r(x_0))}.
		\end{align*}
	For $\sigma r\le|y-x_0|\le r$, applying \eqref{v-h}, \eqref{h_1} and \eqref{h} in the same manner as the previous estimate yields that
		\begin{align*}
			&(\sigma r)^{2s-1}\int_{B_r(x_0)\setminus B_{\sigma r}(x_0)}\frac{|h-\ell_h|}{|x_0-y|^{n+2s}}dy\\
			=& (\sigma r)^{2s-1}\int_{B_{r/2}(x_0)\setminus B_{\sigma r}(x_0)}\frac{|h-\ell_h|}{|x_0-y|^{n+2s}}dy+(\sigma r)^{2s-1}\int_{B_r(x_0)\setminus B_{r/2}(x_0)}\frac{|h-\ell_h|}{|x_0-y|^{n+2s}}dy\\
			\le& C\sigma^{\bar{\alpha}}\left(G(u;x_0,r)+r^{2s-\frac{n}{p}-1}\|f\|_{L^p(B_r(x_0))}\right)+C\sigma ^{2s-1}r^{-1}\left(\fint_{B_r(x_0)}|h|dy+\fint_{B_r(x_0)}|\ell_h|dy\right)\\
			\le& C\left(\sigma^{\bar{\alpha}}+\sigma^{2s-1}\right)\left(G(u;x_0,r)+r^{2s-\frac{n}{p}-1}\|f\|_{L^p(B_r(x_0))}\right).
		\end{align*}Then it holds that
		\begin{align}\label{tail}
			&(\sigma r)^{2s-1}\int_{B_r(x_0)\setminus B_{\sigma r}(x_0)}\frac{|u-l_{x_0,r}-\ell_h|}{|x_0-y|^{n+2s}}dy\nonumber\\
			\le&(\sigma r)^{2s-1}\int_{B_r(x_0)\setminus B_{\sigma r}(x_0)}\frac{|u-l_{x_0,r}-h|}{|x_0-y|^{n+2s}}dy+(\sigma r)^{2s-1}\int_{B_r(x_0)\setminus B_{\sigma r}(x_0)}\frac{|h-\ell_h|}{|x_0-y|^{n+2s}}dy\nonumber\\
			\le& \sigma^{-1}r^{-1}\|u-l_{x_0,r}-h\|_{L^\infty(B_r(x_0))}+C\left(\sigma^{\bar{\alpha}}+\sigma^{2s-1}\right)\left(G(u;x_0,r)+r^{2s-\frac{n}{p}-1}\|f\|_{L^p(B_r(x_0))}\right)\nonumber\\
			\le& C\sigma^{-1}r^{2s-\frac{n}{p}-1}\|f\|_{L^p(B_r(x_0))}+C\left(\sigma^{\bar{\alpha}}+\sigma^{2s-1}\right)G(u;x_0,r).
		\end{align}
Finally, a combination of \eqref{h_2} and \eqref{tail} yields that
	\begin{equation*}
		G(u;x_0,\sigma r)\le \mathcal{E}_{l_{x_0,r}+\ell_h}(u;x_0,\sigma r)\le C\left(\sigma^{\bar{\alpha}}+\sigma^{2s-1}\right) G(u;x_0,r)+C\sigma^{-1} r^{2s-\frac{n}{p}-1}\|f\|_{L^p(B_r(x_0))}.
	\end{equation*}Therefore, we complete the proof of Lemma \ref{Duiter}.
\end{proof}
	Below, we prove the main result concerning the potential estimates of gradients by using the decay estimate we just derived above.
	\begin{proof}[Proof of Theorem \ref{thm1.4}]
		By choosing $\sigma=\sigma(n,s,p,\lambda,\Lambda) \in (0, 1)$ small enough such that $C \left(\sigma^{\bar{\alpha}}+\sigma^{2s-1}\right) \le 1/2$ in Lemma \ref{Duiter}, we sum the iterative excess bounds over all scales $r_k=\sigma^k R$, $k\in\mathbb{N}$, as same as the proof of Theorem \ref{thm1.1} and further obtain
		\begin{equation}\label{eq2}
			\sum_{k=0}^\infty G(u;x_0,r_{k})\le G(u;x_0,R)+C\mathbf{I}^{f}_{2s-1,p}(x_0,R).
		\end{equation}
		Assuming that the optimal affine function which reaches the infimum on each ball $B_{r_k}(x_0)$ is $l_k(x):=a_k+A_k\cdot (x-x_0)$, we know that
		\begin{equation*}
			\begin{aligned}
				|A_{k+1}-A_k|\le \frac{C}{r_{k+1}}\fint_{B_{r_{k+1}}(x_0)}|l_{k+1}-l_k|dy&\le \frac{C}{r_{k+1}}\left(\fint_{B_{r_{k+1}}(x_0)}|l_{k+1}-u|dy+\fint_{B_{r_{k}}(x_0)}|u-l_k|dy\right)  \\
				&\le C\left(G(u;x_0,r_{k+1})+G(u;x_0,r_k)\right)
			\end{aligned}
		\end{equation*}and
		\begin{equation*}
			\begin{aligned}
				|A_0|\le \frac{C}{R}\fint_{B_R(x_0)}|l_0|dy&\le \frac{C}{R}\left(\fint_{B_R(x_0)}|l_0-u|dy+\fint_{B_R(x_0)}|u|dy\right)\\
				&\le CG(u;x_0,R)+\frac{C}{R}\fint_{B_R(x_0)}|u|dy\\
				&\le  C\left(\frac{1}{R}\fint_{B_R(x_0)}|u|dy+\frac{1}{R}\operatorname{Tail}(u;x_0,R)\right),
			\end{aligned}
		\end{equation*}where $C=C(n)>0$.
		Then it follows that
		\begin{equation}\label{A_k}
			\sum_{k=0}^{\infty}|A_{k+1}-A_k|\le C\sum_{k=0}^{\infty}G(u;x_0,r_k)\le CG(u;x_0,R)+C\mathbf{I}^{f}_{2s-1,p}(x_0,R).
		\end{equation} Through boundedness of $G(u;x_0,R)$ and $\mathbf{I}^{f}_{2s-1,p}(x_0,R)$ with $p>n/(2s-1)$, we conclude that $\{A_k\}$ is a Cauchy sequence and its limit existence is denoted as $A_\infty$.
	
	It remains to show that $u$ is differentiable at $x_0$ and $Du(x_0)=A_\infty$. We denote $\tilde{u_k}:=u-l_k$ satisfies $\mathcal{I}(\tilde{u_k},x)=f(x)$ in $B_{r_k}(x_0)$ in the viscosity sense. Applying Theorem \ref{thm1.1} at each $z\in B_{r_k/2}(x_0)$ and combining with the fact that $B_{r_k/2}(z)\subset B_{r_k}(x_0)$ and
	\begin{equation*}
		\mathbf I^f_{2s,p}(z,r_k/2)
		\le Cr_k^{2s-\frac np}\|f\|_{L^p(B_{r_k}(x_0))},
	\end{equation*}we obtain
		\begin{equation*}
			\begin{aligned}
				\|\tilde{u_k}\|_{L^\infty(B_{r_k/2}(x_0))}
				&\le C\left(
				\fint_{B_{r_k}(x_0)}|\tilde{u_k}|dy
				+\operatorname{Tail}(\tilde{u_k};x_0,r_k)
				+r_k^{2s-\frac np}\|f\|_{L^p(B_{r_k}(x_0))}
				\right)\\
				&=C\left(r_kG(u;x_0,r_k)
				+r_k^{2s-\frac np}\|f\|_{L^p(B_{r_k}(x_0))}\right).
			\end{aligned}
		\end{equation*}Let $\varepsilon_k
:=r_k^{-1}\|\tilde{u_k}\|_{L^\infty(B_{r_k/2}(x_0))}.$
The preceding estimate leads to
\begin{equation*}
	\varepsilon_k
	\le CG(u;x_0,r_k)
	+Cr_k^{2s-1-\frac np}\|f\|_{L^p(B_R(x_0))}
	\to 0\quad\text{as }k\to \infty,
\end{equation*}by $\sum\limits_{k=0}^{\infty}G(u;x_0,r_k)<\infty$ and \(p>n/(2s-1)\). In particular, the continuity of $u$ yields that
\begin{equation*}
	|u(x_0)-a_k|\le r_k\varepsilon_k.
\end{equation*}
For any sufficiently small nonzero $|h|$ with $h\in\mathbb R^n$, there exists $k\in\mathbb{N}$ such that $\frac{r_{k+1}}2\le |h|<\frac{r_k}2.$
Then $x_0+h\in B_{r_k/2}(x_0)$ and it follows that
\begin{equation*}
	\begin{aligned}
		\frac{|u(x_0+h)-u(x_0)-A_\infty\cdot h|}{|h|}
		&\le
		\frac{|u(x_0+h)-\ell_k(x_0+h)|+|a_k-u(x_0)|}{|h|}
		+|A_k-A_\infty|\\
		&\le
		\frac4{\sigma}\varepsilon_k+|A_k-A_\infty|.
	\end{aligned}
\end{equation*}Evidently, as $k\to \infty$, we must have $|h|\to 0$, which implies that
\begin{equation*}
	\lim_{|h|\to 0}\frac{|u(x_0+h)-u(x_0)-A_\infty\cdot h|}{|h|}\to 0.
\end{equation*}
Thus $u$ is Fr\'{e}chet differentiable at $x_0$, with $Du(x_0)=A_\infty.$ Finally we obtain
		\begin{equation*}
			|Du(x_0)|=|A_{\infty}|\le |A_0|+\sum_{k=0}^{\infty}|A_{k+1}-A_k|\le \frac{C}{R}\left(\fint_{B_R(x_0)}|u|dy+\operatorname{Tail}(u;x_0,R)\right)+C\mathbf{I}^{f}_{2s-1,p}(x_0,R)
		\end{equation*} and we finish the proof of Theorem \ref{thm1.4}.
	\end{proof}
In the sequel, we present a key ingredient for establishing the universal potential estimate for the gradient of the viscosity solution to equation \eqref{prob1}, namely a pointwise estimate for the following truncated nonlocal fractional sharp maximal function constructed via affine functions
\begin{equation*}
	\mathbf{\tilde{N}}_{t}^{\#}[u](x_0,R):=\sup\limits_{0<r\le R}r^{-t}G(u;x_0,r),
\end{equation*}where the exponent  $0\le t<1$ and the ball $B_R(x_0)\Subset\Omega$.
	\begin{theorem}\label{max}
		Let $s\in (1/2,1)$ and $f\in L^p(\Omega)\cap C(\Omega)$ with $p>n/(2s-1)$. Assume that $u\in C(\Omega)\cap L^1_{2s}(\mathbb{R}^n)$ is a viscosity solution to \eqref{prob1} with the form \eqref{i}, $\mathcal{I}\in \mathfrak{J}^1_s(\lambda,\Lambda)$. For any $\tilde{\alpha}\in(0,\min\{\bar{\alpha},2s-1\})$, there exists a positive constant $C=C(n,s,p,\lambda,\Lambda,\omega(\cdot),\tilde{\alpha})$ such that
		\begin{equation*}
			\mathbf{\tilde{N}}_{t}^{\#}[u](x_0,R)\le C\left(\mathbf{M}_{2s-1-t,p}[f](x_0,R)+R^{-t}G(u;x_0,R)\right),
		\end{equation*}holds for any $0\le t \le\tilde{\alpha}$ and  $B_R(x_0)\Subset\Omega$, where the exponent $\bar{\alpha}=\bar{\alpha}(n,s,\lambda,\Lambda)$ is given by Lemma \ref{1+alpha}.
	\end{theorem}
	\begin{proof}
		According to Lemma \ref{Duiter}, for any $0<\rho\le R$, we have
		\begin{equation*}
			G(u;x_0,\sigma\rho)\le C\left(\sigma^{\bar{\alpha}}+\sigma^{2s-1}\right)G(u;x_0,\rho)+C\sigma^{-1}\rho^{2s-1-\frac{n}{p}}\|f\|_{L^p(B_{\rho}(x_0))}.
		\end{equation*}Set $0\le t \le\tilde{\alpha}<\tau:=\min\{\bar{\alpha},2s-1\}$ and by multiplying both sides by $(\sigma\rho)^{-t}$, there exists constant $C=C(n,s,p,\lambda,\Lambda,\omega(\cdot))>0$ such that
		\begin{equation}\label{eq1}
			(\sigma\rho)^{-t}G(u;x_0,\sigma\rho)\le C\sigma^{\tau-t}\rho^{-t}G(u;x_0,\rho)+C\sigma^{-1-t}\rho^{2s-1-\frac{n}{p}-t}\|f\|_{L^p(B_{\rho}(x_0))}.
		\end{equation}We choose a sufficiently small $\sigma=\sigma(n,s,p,\lambda,\Lambda,\tilde{\alpha})\in(0,1)$ such that $C\sigma^{\tau-t}\le 1/2$. Meanwhile, a direct calculation deduce that
		\begin{equation*}
			\rho^{2s-1-\frac{n}{p}-t}\|f\|_{L^p(B_{\rho}(x_0))}\le C\rho^{2s-1-t}\left(\fint_{B_\rho(x_0)}|f(y)|^pdy\right)^{1/p}\le C\mathbf{M}_{2s-1-t,p}[f](x_0,R),
		\end{equation*}where $C=C(n,p)>0$. Then it follows from \eqref{eq1} that
	\begin{equation*}
		(\sigma\rho)^{-t}G(u;x_0,\sigma\rho)\le \frac{1}{2}\rho^{-t}G(u;x_0,\rho)+C\mathbf{M}_{2s-1-t,p}[f](x_0,R).
	\end{equation*}Setting $r_k:=\sigma^k R$ with $k\in\mathbb{N}$, we obtain the recurrence formula
	\begin{equation*}
		r_{k+1}^{-t}G(u;x_0,r_{k+1})\le \frac{1}{2}r_{k}^{-t}G(u;x_0,r_{k})+C\mathbf{M}_{2s-1-t,p}[f](x_0,R),
	\end{equation*}and further derive
	\begin{equation*}
		\begin{aligned}
			r_{k}^{-t}G(u;x_0,r_{k})&\le \frac{1}{2^k}R^{-t}G(u;x_0,R)+C\mathbf{M}_{2s-1-t,p}[f](x_0,R)\left(\sum_{j=0}^{k-1}\frac{1}{2^j}\right)\\
			&\le R^{-t}G(u;x_0,R)+C\mathbf{M}_{2s-1-t,p}[f](x_0,R).
		\end{aligned}
	\end{equation*}Thus, we conclude that
	\begin{equation*}
		\sup_{k\in\mathbb{N}}r_{k}^{-t}G(u;x_0,r_{k})\le R^{-t}G(u;x_0,R)+C\mathbf{M}_{2s-1-t,p}[f](x_0,R).
	\end{equation*}For any $0<r\le R$, there exists a positive integer $k$ such that $r_{k+1}\le r<r_k$, then we get
	\begin{equation*}
		r^{-t}G(u;x_0,r)\le C r_k^{-t}G(u;x_0,r_k)\le C\sup_{k\in\mathbb{N}}r_k^{-t}G(u;x_0,r_k).
	\end{equation*}Finally, by taking the supremum of the left hand term in above inequality with respect to $r$, we have
	\begin{equation*}
		\sup_{0<r\le R}r^{-t}G(u;x_0,r)\le C\sup_{k\in\mathbb{N}}r_k^{-t}G(u;x_0,r_k)\le C\left(R^{-t}G(u;x_0,R)+\mathbf{M}_{2s-1-t,p}[f](x_0,R)\right),
	\end{equation*}where the positive constant $C$ depends only on $n, s, p, \lambda, \Lambda,\tilde{\alpha}$ and the modulus of continuity $\omega$.
	\end{proof}
Now we can derive an oscillatory estimate of the gradient of the solutions to the nonlocal fully nonlinear elliptic equation \eqref{prob1} with the form \eqref{i}.
	\begin{proof}[Proof of Theorem \ref{thm1.5}]
		For any $x,y\in B_{R/8}(x_0)$ and denote $d:=|x-y|<R/4$. By applying Theorem \ref{max} with the fractional order $0\le \alpha \le \tilde{\alpha}<\tau:=\min\{\bar{\alpha},2s-1\}$, we already have the pointwise bound for the global maximal function of the gradient
		\begin{equation}\label{Mz}
			\mathbf{\tilde{N}}^\#_{\alpha}[u](z, R/2) \le C \left( \mathbf{M}_{2s-1-\alpha, p}[f](z,R/2) + R^{-\alpha} G(u; z, R/2) \right),
		\end{equation}for any point $z \in B_{R/4}(x_0)$. We decompose the gradient difference as follows
		\begin{equation*}
			|Du(x)-Du(y)|\le |Du(x)-A_d^x|+|A_d^x-A_d^y|+|Du(y)-A_d^y|,
		\end{equation*}where $A_d^x$ is the slope of the affine function $l_{x,d}$ corresponding to the infimum of $G(u;x,d)$ and $Du(x)=\lim\limits_{d\to0}A_d^x$. By an argument analogous to that for \eqref{A_k}, we obtain
		\begin{equation*}
			\begin{aligned}
				|Du(x)-A_d^x|&\le \sum_{j=0}^{\infty}|A_{d/2^j}^x-A_{d/2^{j+1}}^x|\\
				&\le CG(u;x,d)+C\mathbf{I}^{f}_{2s-1,p}(x,d)\\
				&\le Cd^{\alpha}\mathbf{\tilde{N}}^\#_{\alpha}[u](x,R/4)+C\mathbf{I}^{f}_{2s-1,p}(x,d).
			\end{aligned}
		\end{equation*}In a similar manner, we deduce that
	\begin{equation}
		|Du(y)-A_d^y|\le Cd^{\alpha}\mathbf{\tilde{N}}^\#_{\alpha}[u](y,R/4)+C\mathbf{I}^{f}_{2s-1,p}(y,d).
	\end{equation}
Since $B_d(x),B_d(y)\subset B_{2d}(x)$, we compute
		\begin{equation*}
			|A_d^x-A_d^y|\le |A_d^x-A_{2d}^x|+|A_{2d}^x-A_d^y|\le CG(u;x,2d)\le Cd^{\alpha}\mathbf{\tilde{N}}^\#_{\alpha}[u](x,R/2).
		\end{equation*}Consequently, by using Lemma \ref{Riesz} and \eqref{Mz}, we derive
			\begin{align}\label{x-y}
				&|Du(x)-Du(y)|\nonumber\\
				\le& Cd^{\alpha}\left(\mathbf{\tilde{N}}^\#_{\alpha}[u](x,R/2)+\mathbf{\tilde{N}}^\#_{\alpha}[u](y,R/2)\right)+C\left(\mathbf{I}^{f}_{2s-1,p}(x,d)+\mathbf{I}^{f}_{2s-1,p}(y,d)\right)\nonumber\\
				\le &Cd^{\alpha}\left(\mathbf{\tilde{N}}^\#_{\alpha}[u](x,R/2)+\mathbf{\tilde{N}}^\#_{\alpha}[u](y,R/2)\right)+Cd^\alpha\left(\mathbf{I}^{f}_{2s-1-\alpha,p}(x,R/2)+\mathbf{I}^{f}_{2s-1-\alpha,p}(y,R/2)\right)\\
				\le& C\left(\frac{d}{R}\right)^\alpha G(u;x_0,R)+Cd^\alpha\left(\mathbf{I}^{f}_{2s-1-\alpha,p}(x,R)+\mathbf{I}^{f}_{2s-1-\alpha,p}(y,R)\right)\nonumber,
			\end{align}where $C=C(n,s,p,\lambda,\Lambda,\omega(\cdot),\tilde{\alpha})>0$. Hence, we finish the proof of Theorem \ref{thm1.5}.
	\end{proof}
Finally, as an application of the gradient potential estimate and the pointwise estimate for the truncated nonlocal fractional sharp maximal function in terms of affine functions, we establish the higher Sobolev regularity estimate of Calder\'{o}n-Zygmund type for viscosity solutions to
nonlocal fully nonlinear ellptic equation \eqref{prob1} with the form \eqref{i}.
	\begin{proof}[Proof of Corollary \ref{cor1.7}]
		Let $\tau=\min\{\bar{\alpha},2s-1\}$, for fixed $t\in[0,\tau)$ and $n/(2s-1)< q<n/(2s-1-t)$, we denote $\tilde{t}=t+\epsilon$, where $\epsilon>0$ is small enough such that $\tilde{t}\in(0,\tau)$. We consider the exponent $$\tilde{q}:=\frac{nq}{n-\beta q}=\frac{nq}{n-(2s-1-\tilde{t})q}<\frac{nq}{n-(2s-1-t)q}=\bar{q},$$where $\beta:=2s-1-\tilde{t}\in (0,n/q)$. Below we further construct an auxiliary exponent $p_0=\frac{1}{2}\left(q+\frac{n}{2s-1}\right)$, which satisfies $p_0\in(n/{(2s-1)},q)$ owing to $q>n/{(2s-1)}$. Since $f\in L^q(B_R(x_0))$, we apply Lemma \ref{L^pmax} to derive
		\begin{equation}\label{eq5}
			\|\mathbf{M}_{\beta,p_0}[f\chi_{B_{R}(x_0)}](\cdot,R/2)\|_{L^{\tilde{q}}(\mathbb{R}^n)}\le C\|f\|_{L^q(B_R(x_0))}.
		\end{equation}
		Then by Theorem \ref{max} and \eqref{eq5}, we obtain
		\begin{equation}\label{M}
				\|\mathbf{\tilde{N}}^{\#}_{\tilde{t}}[u](\cdot,R/2)\|_{L^{\tilde{q}}(B_{R/2}(x_0))}\le C\left(R^{\frac{n}{\tilde{q}}-\tilde{t}}G(u;x_0,R)+\|f\|_{L^{q}(B_{R}(x_0))}\right).
		\end{equation}Meanwhile, according to \cite[Theorem 2.3]{Wolff} and \eqref{eq5}, we deduce that
	\begin{align}\label{I}
		\|\mathbf{I}^f_{\beta,p_0}(\cdot,R/2)\|_{L^{\tilde{q}}(B_{R/2}(x_0))}\le C\|\mathbf{M}_{\beta,p_0}[f\chi_{B_{R}(x_0)}](\cdot,R/2)\|_{L^{\tilde{q}}(\mathbb{R}^n)}\le C\|f\|_{L^q(B_R(x_0))},
	\end{align}where $C=C(n,s,q,t)>0$. By using \eqref{x-y}, for $0<r\le R/16$, we have
			\begin{align}\label{osc}
				&r^{-\tilde{t}}\fint_{B_{r}(x)}|Du(y)-(Du)_{B_r(x)}|dy\nonumber\\
				\le& 2r^{-\tilde{t}}\fint_{B_{r}(x)}|Du(y)-Du(x)|dy\nonumber\\
				\le& C\fint_{B_r(x)}\left(\mathbf{\tilde{N}}^\#_{\tilde{t}}[u](x,R/2)+\mathbf{I}^{f}_{\beta,p_0}(x,R/2)\right)+\left(\mathbf{\tilde{N}}^\#_{\tilde{t}}[u](y,R/2)+\mathbf{I}^{f}_{\beta,p_0}(y,R/2)\right)dy\nonumber\\
				\le& C\left(\mathbf{\tilde{N}}^\#_{\tilde{t}}[u](x,R/2)+\mathbf{I}^{f}_{\beta,p_0}(x,R/2)+\mathbf{M}\left[\mathbf{\tilde{N}}^\#_{\tilde{t}}[u]+\mathbf{I}^{f}_{\beta,p_0}\right](x,R/2)\right).
			\end{align}Given that
	\begin{equation*}
		\beta+\frac{n}{\tilde{q}}-\frac{n}{q}=0\quad\text{and}\quad\tilde{t}-\frac{n}{\tilde{q}}=t-\frac{n}{\bar{q}},
	\end{equation*}based on \eqref{M}-\eqref{osc} and the boundedness of Hardy-Littlewood maximal operator, we conclude that
			\begin{align}\label{eq6}
				&\|\mathbf{M}^{\#}_{\tilde{t}}[Du](\cdot,R/16)\|_{L^{\tilde{q}}(B_{R/8}(x_0))}\nonumber\\
				\le& C\|\mathbf{\tilde{N}}^\#_{\tilde{t}}[u](\cdot,R/2)+\mathbf{I}^{f}_{\beta,p_0}(\cdot,R/2)\|_{L^{\tilde{q}}(B_{R/4}(x_0))}\nonumber\\
				\le& C\left(\|f\|_{L^q(B_R(x_0))}+R^{\frac{n}{\tilde{q}}-\tilde{t}}G(u;x_0,R)\right)\nonumber\\
				\le & CR^{\frac{n}{\bar{q}}-{t}-1}\left(\fint_{B_R(x_0)}|u|dy+\operatorname{Tail}(u;x_0,R)\right)+C\|f\|_{L^{q}(B_{R}(x_0))}.
			\end{align}
	 Furthermore, according to Theorem \ref{thm1.4} with $p=p_0$, we derive
		\begin{equation}\label{u}
			\|Du\|_{L^{\tilde{q}}(B_{R/2}(x_0))}\le CR^{\frac{n}{\tilde{q}}-1}\left(\fint_{B_R(x_0)}|u|dy+\operatorname{Tail}(u;x_0,R)\right)+CR^{\tilde{t}}\|f\|_{L^{q}(B_{R}(x_0))}.
		\end{equation}
		Finally, combining \eqref{eq6} and \eqref{u}, we apply Lemma \ref{Wtq}, thereby obtaining the local Calder\'{o}n-Zygmund estimate as follows
		\begin{equation*}
			\begin{aligned}
				&R^{-t}\|Du\|_{L^{\bar{q}}(B_{R/16}(x_0))}+[Du]_{W^{t,\bar{q}}(B_{R/16}(x_0))}\\
				\le& CR^{-\tilde{t}}\|Du\|_{L^{\tilde{q}}(B_{R/16}(x_0))}+C\|\mathbf{M}^{\#}_{\tilde{t}}[Du](\cdot,R/16)\|_{L^{\tilde{q}}(B_{R/8}(x_0))}\\
				\le& CR^{\frac{n}{\bar{q}}-{t}-1}\left(\fint_{B_R(x_0)}|u|dy+\operatorname{Tail}(u;x_0,R)\right)+C\|f\|_{L^{q}(B_{R}(x_0))},
			\end{aligned}
		\end{equation*}where $C=C(n,s,q,\lambda,\Lambda,\omega(\cdot),t)>0$. Therefore, the proof of Corollary \ref{cor1.7} is complete.
	\end{proof}
	\section{Further results}\label{S6}
	In this section, we particularly consider the concave fully nonlinear integro-differential operator with the following form
	\begin{equation}\label{ii}
		\mathcal{I}(u,x):=\inf_{\beta\in \mathcal{B}} \{-\mathcal{L}_{\beta}u(x)\},
	\end{equation}with $\mathcal{L}_\beta\in\mathfrak{L}^2_s(\lambda,\Lambda)$. Since the proofs of the main results follow the same overall strategy as those developed for nonlocal fully nonlinear equation \eqref{prob1} with the form \eqref{i} in Section \ref{section4}, we omit the repetitive arguments and refer the reader to that section whenever the corresponding steps remain unchanged. We focus instead on the points at which the additional concavity structure of the nonlocal fully nonlinear operator yields an improvement in the admissible regularity range.

Firstly, we present the underlying reason why regularity can be improved under a concave structure and the kernel satisfying $\eta$-H\"{o}lder condition, by virtue of the following $C^{2s,\eta}$ estimates for the solutions.
\begin{lemma}\label{2s+alpha}
	Let $s\in(1/2,1)$ and $0<\lambda\le \Lambda$. There exists $\alpha_1>0$ depending only on $n,s,\lambda,\Lambda$. For some $\eta\in (0,\alpha_1)$ such that $2s+\eta\notin\mathbb{N}$, assume that $u\in C(B_1)\cap L^1_{2s}(\mathbb{R}^n)$ is a viscosity solution of $\mathcal{I}(u,x)=0$ in $B_1$ with the form \eqref{ii} and $\mathcal{I}\in \mathfrak{J}^2_s(\lambda,\Lambda)$. Then $u\in C^{2s,\eta}_{loc}(B_1)$ and
	\begin{equation*}
		\|u\|_{C^{2s,\eta}(B_{1/2})}\le C\|u\|_{L^1_{2s}(\mathbb{R}^n)},
	\end{equation*} where the positive constant $C=C(n,s,\lambda,\Lambda,\eta,M)$.
\end{lemma}
\begin{proof}
	As in the proof of Lemma \ref{1+alpha}, we need a cutoff function $\psi\in C^\infty_c(B_{4/5})$ satisfying $0\le \psi\le 1$ and $\psi\equiv 1$ in $B_{3/4}$. We define $u_1:=\psi u$ and $u_2:=u-u_1$. By using Lemma \ref{sup}, we derive
	\begin{equation}\label{u11}
		\|u_1\|_{L^\infty(\mathbb{R}^n)}\le C\|u\|_{L^1_{2s}(\mathbb{R}^n)}.
	\end{equation}Denote $c_\beta(x):=-\mathcal{L}_{\beta}u_2(x)$ for $\beta\in\mathcal{B}$. Then it follows that
\begin{equation*}
	\tilde{\mathcal{I}}(u_1,x):=\inf_{\beta\in \mathcal{B}}\{-\mathcal{L}_\beta u_1(x)+c_\beta(x)\}=\mathcal{I}(u,x)=0\quad\text{in }B_{2/3}
\end{equation*}in the viscosity sense. Since $u_2=0$ in $B_{3/4}$, upon a change of variables we obtain
\begin{equation*}
	\begin{aligned}
		c_\beta(x)&=\int_{\mathbb{R}^n}u_2(z)\frac{(z-x)^TA_\beta(x)(z-x)}{|z-x|^{n+2s+2}}dz\\
		&=\sum_{i=1}^{n}\sum_{j=1}^{n}A_{ij}(x)T_{ij}(x),
	\end{aligned}
\end{equation*}where $A_\beta(x)=(A_{ij}(x))_{n\times n}$, $T(x)=(T_{ij}(x))_{n\times n}$ and denote $$T_{ij}(x):=\int_{\mathbb{R}^n\setminus B_{3/4}}u_2(z)\frac{(z_i-x_i)(z_j-x_j)}{|z-x|^{n+2s+2}}dz.$$
By direct calculation, it holds that
\begin{equation*}
	\sup_{x\in B_{2/3}}|T(x)|+\sup_{x\in B_{2/3}}|D_xT(x)|\le C\|u_2\|_{L^1_{2s}(\mathbb{R}^n)}\le C\|u\|_{L^1_{2s}(\mathbb{R}^n)}.
\end{equation*}
According to Cauchy-Schwarz inequality, for any $x,x'\in B_{2/3}$, we have
	\begin{equation*}
		\begin{aligned}
				|c_\beta(x)-c_\beta(x')|&\le |A_\beta(x)-A_\beta(x')||T(x)|+|A_\beta(x')||T(x)-T(x')|\\
				&\le C\|u\|_{L^1_{2s}(\mathbb{R}^n)}\left(M|x-x'|^\eta+\Lambda|x-x'|\right).
		\end{aligned}
	\end{equation*}Then it holds that
\begin{equation}\label{c_beta}
	\sup_{\beta\in\mathcal{B}}\|c_\beta\|_{C^\eta(B_{2/3})}\le C(M+\Lambda)\|u\|_{L^1_{2s}(\mathbb{R}^n)}.
\end{equation}Finally, based on \eqref{u11} and \eqref{c_beta}, we use \cite[Theorem 1.3(b)]{serra} with the viscosity solution $u_1$ of $\tilde{\mathcal{I}}(u_1,x)=0$ in $B_{2/3}$ to obtain that
\begin{equation*}
	\|u\|_{C^{2s,\eta}(B_{1/2})}=\|u_1\|_{C^{2s,\eta}(B_{1/2})}\le C\left(\|u_1\|_{L^\infty(\mathbb{R}^n)}+\|u\|_{L^1_{2s}(\mathbb{R}^n)}\right)\le C\|u\|_{L^1_{2s}(\mathbb{R}^n)},
\end{equation*}where the positive constant $C$ depends on $n,s,\lambda,\Lambda,\eta$ and $M$.
\end{proof}
The concavity structure leads to improved regularity of the solution, based on Lemma \ref{2s+alpha}, a proof similar to Lemma \ref{Duiter} yields the refined decay estimate for the global excess functional constructed by affine functions.
\begin{lemma}
	Let $s\in(1/2,1)$, $\sigma\in(0,1)$ and  $f\in L^p(\Omega)\cap C(\Omega)$, $p>n/(2s-1)$. Assume that $u\in C(\Omega)\cap L^1_{2s}(\mathbb{R}^n)$ is a viscosity solution to $\mathcal{I}(u,x)=f$ in $\Omega$ of the form \eqref{ii} with $\mathcal{I}\in \mathfrak{J}^2_s(\lambda,\Lambda)$. Then there exists constant $C=C(n,s,p,\lambda,\Lambda,\eta,M)>0$ such that
	\begin{equation*}
		G(u;x_0,\sigma R)\le C\sigma^{2s-1}G(u;x_0,R)+C\sigma^{-1} R^{2s-\frac{n}{p}-1}\|f\|_{L^p(B_R)}
	\end{equation*}holds for $B_{\sigma R}(x_0)\subset B_{R}(x_0)\Subset \Omega$, where $\eta\in (0,\alpha_1)$ is given by Lemma \ref{2s+alpha}.
\end{lemma}
\begin{proof}
	In the proof of Lemma \ref{Duiter}, we generalize the solution of the homogeneous equation. For the concave form \eqref{ii}, based on Lemma \ref{2s+alpha}, we obtain a higher regularity estimate for $h$
	\begin{equation*}
		\|h\|_{C^{2s,\eta}(B_{R/2})}\le CR^{-2s-\eta}\|h_R\|_{L^1_{2s}(\mathbb{R}^n)},
	\end{equation*}where $\eta\in (0,\alpha_1)$ such that $2s+\eta\notin\mathbb{N}$. Finally, we repeat the above proof process of Lemma \ref{Duiter} to obtain
	\begin{equation}\label{eq7}
		G(u;x_0,\sigma R)\le  C\sigma^{2s-1} G(u;x_0,R)+C\sigma^{-1} R^{2s-\frac{n}{p}-1}\|f\|_{L^p(B_R)},
	\end{equation}where $C=C(n,s,p,\lambda,\Lambda,\eta,M)>0$.
\end{proof}
Next we establish the pointwise estimate for fractional sharp maximal function of the global affine excess functional. From \eqref{eq7}, we directly obtain the following property by referring to the proof of Theorem \ref{max}.
\begin{theorem}
	Let $s\in (1/2,1)$ and $f\in L^p(\Omega)\cap C(\Omega)$ with $p>n/(2s-1)$. Assume that $u\in C(\Omega)\cap L^1_{2s}(\mathbb{R}^n)$ is a viscosity solution to \eqref{prob1} with the form \eqref{ii}, $\mathcal{I}\in \mathfrak{J}^2_s(\lambda,\Lambda)$. For any $\tilde{\alpha}\in(0,2s-1)$, there exists a positive constant $C=C(n,s,p,\lambda,\Lambda,\eta,M,\tilde{\alpha})$ such that
	\begin{equation*}
		\mathbf{\tilde{N}}_{t}^{\#}[u](x_0,R)\le C\left(\mathbf{M}_{2s-1-t,p}[f](x_0,R)+R^{-t}G(u;x_0,R)\right)
	\end{equation*}holds for any $0\le t\le\tilde{\alpha}$ and  $B_R(x_0)\Subset\Omega$, where the  exponent $\eta\in(0,\alpha_1)$ is given by Lemma \ref{2s+alpha}.
\end{theorem}
Once the above improved estimates are obtained, the oscillatory estimates concerning the gradient of the viscosity solution follow from the same iteration, telescoping argument and harmonic analysis tools used in the proof of Theorem \ref{thm1.5}.
\begin{theorem}\label{thm1.6}
	Let $s\in (1/2,1)$ and $f\in L^p(\Omega)\cap C(\Omega)$ with $p>n/(2s-1)$. Assume that $u\in C(\Omega)\cap L^1_{2s}(\mathbb{R}^n)$ is a viscosity solution to \eqref{prob1} with the form \eqref{ii}, $\mathcal{I}\in \mathfrak{J}^2_s(\lambda,\Lambda)$. Then for any $\tilde{\alpha}\in(0,2s-1)$, there exists a constant $C=C(n,s,p,\lambda,\Lambda,\eta,M,\tilde{\alpha})>0$ such that
	\begin{equation*}
		\begin{aligned}
				|Du(x)&-Du(y)|\\&\le C\left[\mathbf{I}^{f}_{2s-1-\alpha,p}(x,R)+\mathbf{I}^{f}_{2s-1-\alpha,p}(y,R)+\frac{1}{R^{1+\alpha}}\left(\fint_{B_{R}(x_0)}|u|dy+\operatorname{Tail}(u;x_0,R)\right)\right]|x-y|^\alpha
		\end{aligned}
	\end{equation*}holds for every $x,y\in B_{R/8}(x_0)\Subset\Omega$ and for any  $0\le\alpha\le\tilde{\alpha}$, where the  exponent $\eta\in(0,\alpha_1)$ is given by Lemma \ref{2s+alpha}.
\end{theorem}
Finally, we raise the integrability exponent $t$ to any value between $0$ and $2s-1$ for the viscosity solutions to the nonlocal fully nonlinear equation \eqref{prob1} with the form \eqref{ii}.
\begin{corollary}
	Let $s\in (1/2,1)$. For any $t\in(0,2s-1)$ assume that $f\in L^q(\Omega)$, $n/(2s-1)< q<n/(2s-1-t)$ and denote $\bar{q}:=\frac{nq}{n-(2s-1-t)q}$. Assume that $u\in C(\Omega)\cap L^1_{2s}(\mathbb{R}^n)$ is a viscosity solution to \eqref{prob1} with the form \eqref{ii}, $\mathcal{I}\in \mathfrak{J}^2_s(\lambda,\Lambda)$. Then $u\in W^{1+t,\bar{q}}_{loc}(\Omega)$. Moreover, for $B_R(x_0)\Subset\Omega$, there exists a positive constant $C=C(n,s,q,\lambda,\Lambda,\eta,M,t)>0$, it holds that
	\begin{align*}
		&R^{-t}\|Du\|_{L^{\bar{q}}(B_{R/16}(x_0))}+[Du]_{W^{t,\bar{q}}(B_{R/16}(x_0))}\\
		\le& CR^{\frac{n}{\bar{q}}-t-1}\left(\fint_{B_R(x_0)}|u|dy+\operatorname{Tail}(u;x_0,R)\right)+C\|f\|_{L^{q}(B_{R}(x_0))}.
	\end{align*}
\end{corollary}

	\section*{Acknowledgments}
	The authors are supported by the National Natural Science Foundation of China (NNSF  Grant No.12571103), Young Scientific and Technological Talents (Level Three) in Tianjin, and Natural Science Foundation of Tianjin (No.25JCQNJC01400).

\end{document}